\documentclass[12pt]{article}

\usepackage{amsmath}
\usepackage{amssymb}
\usepackage{amsfonts}
\usepackage{a4}
\usepackage{graphicx}

\newtheorem{theorem}{Theorem}[section]
\newtheorem{lemma}[theorem]{Lemma}

\newtheorem{definition}[theorem]{Definition}

\newcommand{\Proof}{\par\noindent{\em Proof. }}
\newcommand{\eop}{\nopagebreak\hspace*{\fill}$\Box$\smallskip}

\newcommand{\Z}{\Bbb Z}
\newcommand{\R}{\Bbb R}

\def\Id{\mathbf{Id}}
\def\id{\mathbf{id}}
\def\calL{\mathcal{L}}
\def\eps{\varepsilon}
\def\vv{\mathbf{v}}
\def\trace{\operatorname{trace}}

\def\e{\mathbf{e}}

\def\weakly{\rightharpoonup}

\def\dist{\operatorname{dist}}

\def\XXint#1#2#3{{\setbox0=\hbox{$#1{#2#3}{\int}$}
     \vcenter{\hbox{$#2#3$}}\kern-.5\wd0}}

\usepackage{color}

\definecolor{dark_green}{rgb}{0, .5, 0}

\begin{document}

\begin{center}
\begin{Large}
{\bf { On a discrete-to-continuum convergence result for a two dimensional brittle material in the small displacement regime}}
\end{Large}
\end{center}

\begin{center}
\begin{large}
Manuel Friedrich\footnote{ Universit{\"a}t Augsburg, Institut f{\"u}r Mathematik, 
Universit{\"a}tsstr.\ 14, 86159 Augsburg, Germany. {\tt manuel.friedrich@math.uni-augsburg.de}}
and Bernd Schmidt\footnote{Universit{\"a}t Augsburg, Institut f{\"u}r Mathematik, 
Universit{\"a}tsstr.\ 14, 86159 Augsburg, Germany. {\tt bernd.schmidt@math.uni-augsburg.de}}\\
\end{large}
\end{center}

\begin{center}
\today
\end{center}
\bigskip

\begin{abstract}

We consider a two-dimensional atomic mass spring system and show that in the small displacement regime the corresponding discrete energies can be related to a continuum Griffith energy functional in the sense of $\Gamma$-convergence. We also analyze the continuum problem for a rectangular bar under tensile boundary conditions and find that depending on the boundary loading the minimizers are either homogeneous elastic deformations or configurations that are completely cracked generically along a crystallographic line. As applications we discuss cleavage properties of strained crystals and an effective continuum fracture energy for magnets. 

\end{abstract}
\bigskip

\begin{small}

\end{small}

\tableofcontents

\section{Introduction}

A fundamental problem in static fracture mechanics is to determine the behavior of a brittle material which is subject to certain displacements imposed at its boundary. Of particular interest is the identification of critical loads at which failure occurs. A natural framework to treat such \textit{free discontinuity problems} with variational methods is given by Griffith energy functionals introduced by Francfort and Marigo \cite{Francfort-Marigo:1998} comprising elastic bulk contributions and surface terms comparable to the size of the crack (see also \cite{DeGiorgi-Ambrosio:1988}). Often these models contain anisotropic surface terms (see e.g.\ \cite{Alicandro-Focardi-Gelli:2000, Focardi-Gelli:2003, Negri:2003}) modeling the fact that due to the crystalline structure of the materials certain directions for the formation of cracks are energetically favored. Indeed, fracture typically occurs in the form of cleavage along crystallographic planes. Ultimately, such a continuum model should be identified as an effective theory derived from atomistic interactions.  

Specifying the set-up even further, a basic experiment to infer material properties of brittle materials is to probe the specimen by applying a uniaxial tensile strain which allows to determine its Poisson ratio in the elastic regime and a critical load beyond which the body fails due to fracture. From a theoretical point of view this problem has been studied recently by Mora-Corral in \cite{Mora:2010}, where he investigates a rectangular bar of brittle, incompressible, homogeneous and isotropic material subject to uniaxial extension and shows that, depending on the loading, the minimizers are either given by purely elastic configurations or deformations with horizontal fracture. 

An atomistic model problem with surface contributions sensitive to the crack geometry has been studied by the authors in \cite{FriedrichSchmidt:2014} leading to a complete analysis of the asymptotically optimal configurations under uniaxial extension in the discrete-to-continuum limit: The body shows pure elastic behavior in the subcritical case and for supercritical boundary values generically cleavage occurs along a specific crystallographic line. However, for a certain symmetric orientation of the lattice cleavage may fail more complicated crack geometries are possible.

The goal of this work is to show that in the small displacement regime the energies associated to such a discrete system can be related to a continuum Griffith energy functional with anisotropic surface contributions in the sense of $\Gamma$-convergence. Moreover, we analyze the continuum problem under tensile boundary conditions. In this way we (1) obtain convergence scheme which in certain applications to be discussed below allows to identify effective continuum fracture energies, (2) extend the results of \cite{Mora:2010} to anisotropic and compressible materials and (3) re-derive in part the aforementioned convergence results of \cite{FriedrichSchmidt:2014}. 

In the theory of fracture mechanics the passage from discrete systems to continuum models via $\Gamma$-convergence is by now well understood for one-dimensional chains, see e.g. \cite{Braides-Cicalese:2007,Braides-DalMaso-Garroni:1999,Braides-Gelli:2002-1}. In the higher dimensional setting there are results for scalar valued models  (see \cite{Braides-Gelli:2002-2}) and approximations of vector valued free discontinuity problems where the elastic bulk part of the energy is characterized by linearized terms (see \cite{Alicandro-Focardi-Gelli:2000}) or by a quasiconvex stored energy density (see \cite{Focardi-Gelli:2003}). However, in more than one dimension the energy density of discrete systems such as well-known mass spring models is in general not given in terms of a discretized continuum quasiconvex function. For large strains these lattices typically become even unstable, see e.g.\ the basic model discussed in \cite{Friesecke-Theil:02}. Consequently, in the regime of finite elasticity it is a subtle question if minimizers for given boundary data exist at all. On the other hand, for sufficiently small strains one may expect the Cauchy-Born rule to apply so that individual atoms do in fact follow a macroscopic deformation gradient, see \cite{Friesecke-Theil:02, Conti-Dolzmann-Kirchheim-Mueller:06}. In particular this applies to the regime of infinitesimal elastic strains. For purely elastic interactions this relation has also been obtained in the sense of $\Gamma$-convergence for a simultaneous passage from discrete to continuum and linearization process in \cite{Braides-Solci-Vitali:07, Schmidt:2009}. 

The model considered in \cite{FriedrichSchmidt:2014} as well as the one-dimesional seminal paper \cite{Braides-Lew-Ortiz:06}  suggest that the most interesting regime for the elastic strains is given by $\sqrt{\eps}$ ($\eps$ denotes the typical interatomic distance) as in this particular regime the elastic and the crack energy are of the same order. This is in accordance to the observation that brittle materials develop cracks already at moderately large strains. Moreover, it shows that a discrete-to-continuum $\Gamma$-limit for the discrete energies under consideration naturally involves a linearization process. 

Identifying all possible limiting continuum configurations and energies is a challenging task as necessary smallness assumptions on the discrete gradient can not be inferred from suitable energy bounds and deriving rigidity estimates being essential in the passage from nonlinear to linearized theory (see \cite{Braides-Solci-Vitali:07, Schmidt:2009}) is a subtle problem. Partial results have been obtain in \cite{FriedrichSchmidt:2014} for almost minimizers of a boundary value problem describing uniaxial extension. A general analysis in two dimensions is deferred to a subsequent work. In the present context we make the simplifying assumption that we consider deformations lying $\sqrt{\eps}$-close to the identity mapping. However, we will also see that there are physically interesting applications e.g.\ to magnetic materials where such an assumption can be justified rigorously. 
 
It then turns out that the derivation of the continuum limit is an issue similar to those considered in \cite{Alicandro-Focardi-Gelli:2000, Braides-Gelli:2002-2, Focardi-Gelli:2003}. Nevertheless, we believe that the present $\Gamma$-convergence result is interesting as (1) it gives rise to a limiting Griffith functional in the realm of linearized elasticity which can be explicitly investigated for cleavage, (2) there are applications   to systems with small displacements for small energies and (3) to the best of our knowledge our approach to the problem differs from techniques which are predominantly used when treating discrete systems in the framework of fracture mechanics.     

The reduction to one-dimensional sections using slicing properties for \textit{(special) functions of bounded variation} turned out to be a useful tool not only to derive general properties of these function spaces but also to study discrete systems and variational approximation of free discontinuity problems. E.g., the original proofs of the main compactness and closure theorems in $SBV$ (see \cite{Ambrosio:89}) as well as the $\Gamma$-convergence results in \cite{Braides-Gelli:2002-2, Focardi-Gelli:2003} make use of this integralgeometric approach. Similar to the fact that there are simplified proofs of these compactness theorems being derived without  the slicing technique (see \cite{Alberti-Mantegazza}), we show that in our framework the lower bound of the $\Gamma$-limit can be achieved in a different way. In fact, we carefully construct the crack shapes of discrete configurations in an explicit way which allows us to directly appeal to lower semicontinuity results for $SBV$ functions.   

The paper is organized as follows. We first introduce our discrete model and state our main results in Section \ref{sec:model-and-main-results}. Here we also briefly discuss how these results shed new light on our findings in \cite{FriedrichSchmidt:2014} on crystal cleavage and study an application to fractured magnets in an external field. 
 
Section \ref{sec:convergence-variational-problem} is devoted to the derivation of the continuum energy functional via $\Gamma$-convergence. The main idea for the lower bound relies on a separation of the energy into elastic and surface contributions by introducing an interpolation with discontinuities on triangles where large expansion occurs. By constructing the set of discontinuity points in a suitable way the surface energy can be estimated using lower semicontinuity results for $SBV$ functions. The elastic part can be treated similarly as in \cite{FrieseckeJamesMueller:02, Schmidt:2009}. 
 
Finally, in Section \ref{sec:limiting-variational-problem} we analyze the continuum problem under tensile boundary values and extend the results obtained in \cite{Mora:2010} to anisotropic and compressible materials. A careful analysis of the anisotropic surface contribution shows that in the generic case there is a unique optimal direction for the formation of fracture, while in a symmetrically degenerate case cleavage fails and all energetically optimal crack geometries can be characterized by specific Lipschitz curves. As in \cite{Mora:2010} the proof makes use of  a qualitative rigidity result for $SBV$ functions (see  \cite{Chambolle-Giacomini-Ponsiglione:2007}) and of the structure theorem on the boundary of sets of finite perimeter by Federer \cite{Federer:1969}.

\section{The model, main results and applications}\label{sec:model-and-main-results}

\subsection{The discrete model}

 Let ${\cal L}$ denote the rotated triangular lattice 
$$ {\cal L} 
   = R_{\cal L} \begin{pmatrix} 1 & \frac{1}{2} \\ 0 & \frac{\sqrt{3}}{2} \end{pmatrix} \Z^2 
   = \{ \lambda_1 \vv_1 + \lambda_2 \vv_2 : \lambda_1, \lambda_2 \in \Z \}, $$ 
where $R_{\cal L} = \begin{footnotesize} \begin{pmatrix} \cos \phi & - \sin \phi \\ \sin \phi & \cos \phi \end{pmatrix} \in SO(2) \end{footnotesize}$  is some rotation and $\vv_1$, $\vv_2$ are the lattice vectors $\vv_1 = R_{\cal L} \mathbf{e}_1$ and $\vv_2 = R_{\cal L}(\frac{1}{2} \mathbf{e}_1 + \frac{\sqrt{3}}{2} \mathbf{e}_2)$, respectively. Without loss of generality we can assume $\phi \in [0, \frac{\pi}{3})$. We collect the basic lattice vectors in the set ${\cal V} = \left\{\vv_1,\vv_2,\vv_2 - \vv_1\right\}$. The macroscopic region $\Omega \subset \R^2$ occupied by the body is supposed to be a bounded domain with Lipschitz boundary. In the reference configuration the positions of the specimen's atoms are given by the points of the scaled lattice $\eps\calL$ that lie within $\Omega$. Here $\eps$ is a small parameter defining the length scale of the typical interatomic distances. 

The deformations of our system are mappings $y : \eps \calL \cap \Omega \to \R^2$. The energy associated to such a deformation $y$ is assumed to be given by nearest neighbor interactions as 
\begin{align}\label{eq:Energy}
  E_{\eps}(y) 
  = \frac{1}{2} \sum_{x,x' \in \eps {\cal L} \cap \Omega \atop |x-x'| = \eps} W \left( \frac{|y(x) - y(x')|}{\eps} \right).  
\end{align} 
Note that the scaling factor $\frac{1}{\eps}$ in the argument of $W$ takes account of the scaling of the interatomic distances with $\eps$. The pair interaction potential $W:[0,\infty) \to [0, \infty]$ is supposed to be of `Lennard-Jones-type': 
\begin{itemize}\label{W-assumptions}
\item[(i)] $W \ge 0$ and $W(r) = 0$ if and only if $r = 1$. 
\item[(ii)] $W$ is continuous on $[0, \infty)$ and $C^2$ in a neighborhood of $1$ with $\alpha := W''(1) > 0$. 
\item[(iii)] $\lim_{r \to \infty} W(r) = \beta > 0$. 
\end{itemize} 

In order to analyze the passage to the limit as $\eps \to 0$ it will be useful to interpolate and rewrite the energy as an integral functional. Let ${\cal C}_{\eps}$ be the set of equilateral triangles $\triangle \subset \Omega$ of sidelength $\eps$ with vertices in $\eps {\cal L}$ and define $\Omega_{\eps} = \bigcup_{\triangle \in {\cal C}_{\eps}} \triangle$. By $\tilde{y} : \Omega_{\eps} \to \R^2$ we denote the interpolation of $y$, which is affine on each $\triangle \in {\cal C}_{\eps}$. The derivative of $\tilde{y}$ is denoted by $\nabla \tilde{y}$, whereas we write $(y)_{\triangle}$ for the (constant) value of the derivative on a triangle $\triangle \in {\cal C}_{\eps}$. Then \eqref{eq:Energy} can be rewritten as 
\begin{align}\label{eq:E-integral}
\begin{split}
  E_{\eps}(y) 
  &= \sum_{\triangle \in {\cal C}_{\eps}} W_{\triangle} ((\tilde{y})_{\triangle}) 
     + E_{\eps}^{\rm boundary}(y) \\ 
  &= \frac{4}{\sqrt{3}\eps^2} \int_{\Omega_{\eps}} W_{\triangle} (\nabla \tilde{y}) \, dx 
     + E_{\eps}^{\rm boundary}(y), 
\end{split} \end{align}
where 
\begin{align}\label{eq:W-triangle}
  W_{\triangle}(F) 
  &= \frac{1}{2} \Big( W(|F \vv_1|) + W(|F \vv_2|) + W(|F (\vv_2 - \vv_1)|) \Big).     
\end{align} 
Here we used that $|\triangle| = \sqrt{3}\eps^2/4$. The boundary term is the sum of pair interaction energies $\frac{1}{4} W(\frac{|y(x) - y(x')|}{\eps})$ or $\frac{1}{2} W(\frac{|y(x) - y(x')|}{\eps})$ over nearest neighbor pairs which form the side of only one or no triangle in ${\cal C}_{\eps}$, respectively. 

Due to the discreteness of the underlying atomic lattice,  Dirichlet boundary conditions have to be imposed in a small neighborhood of the boundary as otherwise cracks near the boundary may become energetically more favorable. Assume that $\tilde{\Omega} \supset \Omega$ is a bounded, open domain in $\R^2$ with Lipschitz boundary defining the Dirichlet boundary $\partial_D\Omega = \partial \Omega \cap \tilde{\Omega}$ of $\Omega$. For $g \in W^{1, \infty}(\tilde{\Omega})$ we define the class of discrete displacements assuming the boundary value $g$ on $\partial_D\Omega$ as 
\begin{align}\label{eq: Ag} 
  {\cal A}_g 
  = \big\{ u : \eps {\cal L} \cap \tilde{\Omega} \to \R^2 : 
    u(x) = g(x) \text{ for } x \in \eps {\cal L} \cap \Omega_{D,\eps} \big\}, 
\end{align}
where $\Omega_{D,\eps} := \lbrace x\in \tilde{\Omega}: \dist(x, \tilde{\Omega}\setminus \Omega) \le \eps \rbrace$. For the corresponding deformations $y = \id + u$ this amounts to requiring $y(x) = x + g(x)$ for $x \in \eps {\cal L} \cap \Omega_{D,\eps}$. Similar as before, we let $\tilde{\cal C}_{\eps}$ be the set of equilateral triangles $\triangle \subset \tilde{\Omega}$ with vertices in $\eps {\cal L}$ and define $\tilde{\Omega}_{\eps} = \bigcup_{\triangle \in \tilde{\cal C}_{\eps}} \triangle$. By $\tilde{y} : \tilde{\Omega}_{\eps} \to \R^2$ we again denote the piecewise affine interpolation of $y$.

It is easy to see that the formation of a crack of finite length resulting from a number of largely deformed triangles scaling with $\frac{1}{\eps}$ leads to an energy contribution to $E_{\eps}$ scaling with $\eps$. The most interesting regime is when the elastic energy contributions to $E_{\eps}$ and the energy cost of a cracked configurations are of the same order. We are thus particularly interested in boundary displacements $g_{\eps}$ scaling with $\sqrt{\eps}$. For then there are also completely elastic deformations for which $E_\eps$ scales with $\eps$, e.g. $E_{\eps}(\id + g_{\eps}) = O(\eps)$.

In order to obtain finite energies and displacements in the limit $\eps \to 0$, we accordingly rescale the displacement field to $u = \frac{1}{\sqrt{\eps}}(y -\id)$ and the energy $E_{\eps}$ to 
$$ {\cal E}_{\eps}(u) := \eps E_{\eps}(y) = \eps E_{\eps}(\id + \sqrt{\eps} u).$$
Moreover, we will assume $u \in {\cal A}_{g_\eps}$ for some $g_\eps \in W^{1,\infty}(\tilde{\Omega})$. 
 
We also introduce the functionals ${\cal E}^\chi_\eps$ which arise from ${\cal E}_\eps$ by replacing $W_\Delta$ by $W_{\Delta, \chi} = W_\Delta + \chi$, where $\chi : \R^{2 \times 2} \to [0,  \infty]$ is a frame indifferent penalty term with $\chi \ge c_\chi >0$ in a neighborhood of $O(2)\setminus SO(2)$ and $ \chi \equiv 0$ in a neighborhood of $SO(2) \cup \{\infty\}$. This term is a mild extra assumption to assure that the orientation of the triangles is preserved in the elastic regime and unphysical effects are avoided.

\subsection{Convergence of the variational problems}\label{sec: con}

Our convergence analysis applies to discrete deformations which may elongate a number scaling with $\frac{1}{\eps}$ of springs very largely, leading to cracks of finite length in the continuum limit. On triangles not adjacent to such essentially broken springs, the defomations are $\sqrt{\eps}$-close to the identity mapping, so that the accordingly rescaled displacements are of bounded $L^2$-norm. Note that the first of these assumptions can be inferred from suitable energy bounds. By way of example, however, we will see that this cannot be true for the displacement estimates in the bulk: The sequence of functionals $({\cal E}_{\eps})_{\eps}$ is not equicoercive. Nevertheless, it is interesting to investigate this regime in order to identify a corresponding continuum functional which describes the system in the realm of Griffith models with linearized elasticity. In fact, below we will discuss two specific models where external fields or boundary conditions break the rotational symmetry whence the sequence $({\cal E}^{\chi}_{\eps})_{ \eps}$ satisfies suitable equicoercivity conditions. 

Recall that the space $SBV(\Omega; \R^2)$, abbreviated as $SBV(\Omega)$ hereafter,  of special functions of bounded variation consists of functions $u \in L^1(\Omega; \R^2)$ whose distributional derivative $Du$ is a finite Radon measure, which splits into an absolutely continuous part with density $\nabla u$ with respect to Lebesgue measure and a singular part $D^j u$ whose Cantor part vanishes and thus is of the form 
$$ D^j u = [u] \otimes \nu_u {\cal H}^1 \lfloor J_u, $$
where ${\cal H}^1$ denotes the one-dimensional Hausdorff measure, $J_u$ (the `crack path') is an ${\cal H}^1$-rectifiable set in $\Omega$, $\nu_u$ is a normal of $J_u$ and $[u] = u^+ - u^-$ (the `crack opening') with $u^{\pm}$ being the one-sided limits of $u$ at $J_u$. If in addition $\nabla u \in L^2(\Omega)$ and ${\cal H}^1(J_u) < \infty$, we write $u \in SBV^2(\Omega)$. See \cite{Ambrosio-Fusco-Pallara:2000} for the basic properties of these function spaces. 

The sense in which discrete displacements are considered convergent to a limiting displacement in SBV is made precise in the following definition. 
\begin{definition}\label{def: convergence} 
Suppose $ u_{\eps} : \eps {\cal L} \cap \tilde{\Omega} \to \R^2$ is a sequence of discrete displacements. We say that $u_{\eps}$ converges to some $u \in SBV^2(\tilde{\Omega})$: $u_{\eps} \to u$, if 
\begin{itemize}
\item[(i)] $\chi_{\tilde{\Omega}_{\eps}} \tilde{u}_{\eps} \to u$ in $L^{1}(\tilde{\Omega})$ 
\end{itemize}
and there exists a sequence ${\cal C}^*_{\eps} \subset \tilde{\cal C}_{\eps}$ with $\# {\cal C}^*_{\eps} \le \frac{C}{\eps}$ for a constant $C$ independent of $\eps$ such that 
\begin{itemize}
\item[(ii)] $\| \nabla \tilde{u}_{\eps} \|_{L^2(\tilde{\Omega} \setminus \cup_{\triangle \in {\cal C}^*_{\eps}} \triangle)} \le C$. 
\end{itemize}
\end{definition}

The main idea will be to separate the energy into elastic and crack surface contributions by introducing a threshold such that triangles $\triangle$ with $(y)_{\triangle}$ beyond that threshold are considered as cracked and $\tilde{y}$ is modified there to a discontinuous function. The treatment of the elastic part draws ideas from \cite{Schmidt:2009} and \cite{FrieseckeJamesMueller:02}. To derive the crack energy, one could use a slicing technique, see, e.g., \cite{Braides-Gelli:2002-2}. Although also possible in our framework, we follow a different approach here: We carefully construct crack shapes of discrete configurations in an explicit way which allows us to directly appeal to lower semicontinuity results for $SBV$ functions in order to derive the main energy estimates. 

Consider the limiting functional 
\begin{equation*}
{\cal E}(u) = 
                    \frac{4}{\sqrt{3}} \int_{\Omega}  \frac{1}{2} Q(e(u)) \, dx  + \int_{J_{u}} \sum_{\vv \in {\cal V}} \frac{2\beta}{\sqrt{3}} |\vv \cdot \nu_{u}|\, d{\cal H}^1
\end{equation*}
for $u \in SBV^2(\tilde{\Omega})$, where $e(u)=\frac{1}{2}\left(\nabla u^{T} + \nabla u\right)$ denotes the symmetric part of the gradient. $Q$ is the linearization of $W_{\triangle}$ about the identity matrix $\Id$ (see Lemma \ref{lemma:W-linearized-properties} for its explicit form). Note that for a displacement field $u$, which is the limit of a sequence $(u_{\eps}) \subset  {\cal A}_{g_{\eps}}$ converging in the sense of Definition \ref{def: convergence}, we get $u = g$ on $\tilde{\Omega} \setminus \Omega$, where $g = L^1\text{-}\lim_{\eps \to 0} g_\eps$. Therefore, if $u|_\Omega$ does not attain the boundary condition $g$ on the Dirichlet boundary $\partial_D \Omega$ (in the sense of traces), this will be penalized in the energy ${\cal E}(u)$. In Section \ref{sec:convergence-variational-problem} we prove the following $\Gamma$-convergence result (see \cite{DalMaso:93} for an exhaustive treatment of $\Gamma$-convergence):

\begin{theorem}\label{th: Gamma-convergence}
\begin{itemize}
\item[(i)] Let $(g_\eps)_\eps \subset W^{1,\infty}(\tilde{\Omega})$ with $\sup_\eps \Vert g_\eps \Vert_{W^{1,\infty}(\tilde{\Omega})} < + \infty$. If $(u_{\eps})_\eps$ is a sequence of discrete displacements with $u_{\eps} \in {\cal A}_{g_\eps}$ and $u_\eps \to u \in SBV^2(\tilde{\Omega})$, then 
$$ \liminf_{\eps \to 0} {\cal E}_{\eps}(u_{\eps}) 
   \ge {\cal E}(u). $$ 

\item[(ii)] For every $u \in SBV^2(\tilde{\Omega})$ and $g \in W^{1,\infty}(\tilde{\Omega})$ with $u=g$ on $\tilde{\Omega} \setminus \Omega$ there is a sequence $(u_{\eps})_\eps$ of discrete displacements such that $u_{\eps} \in {\cal A}_{g}$, $u_{\eps} \to u \in SBV^2(\tilde{\Omega})$ and 
$$ \lim_{\eps \to 0} {\cal E}^{ \chi}_{\eps}(u_{\eps}) 
   = {\cal E}(u). $$ 
\end{itemize}
\end{theorem} 
Note that the recovery sequence is obtained for the energy ${\cal E}^\chi_\eps$ which includes the frame indifferent penalty term.
Due to the frame indifference of $W$, $({\cal E}_{\eps})$ and $({\cal E}^\chi_{\eps})$ are not equicoercive as the following example shows. 
\smallskip 

\noindent{\bf Example.} 
Assume that the specimen satisfying the boundary conditions is broken into three parts by two even cracks where the middle part is subject to a rotation $R \neq \Id$ so that 
$$ \nabla \tilde{y}_{\eps} = R \text{ for } p \leq x_1 \leq q, \ 0 < p < q < l. $$
In particular, the energy of the configuration is of order $1$. But for $p \leq x_1 \leq q$
$$ |\nabla \tilde{u}_{\eps}(x)| 
   = \left|\frac{1}{\sqrt{\eps}} \left(R - \Id\right)\right| 
   \rightarrow \infty \text{ for } \eps \rightarrow 0.$$ 
Thus, $\nabla \tilde{u}_{\eps}$ is not bounded in $L^1$ and so $u_{\eps}$ does not converge.
\smallskip

We now add a term to ${\cal E}_{\eps}$ such that the sequence becomes equicoercive. Let $\hat{m}: \R^{2 \times 2} \to S^1$ be a function satisfying 
$$\hat{m}(RF) = R\hat{m}(F) \ \text{ for all } F \in \R^{2 \times 2}, R \in SO(2), \ \ \ \ \hat{m}(\Id) = \e_1 .$$ Moreover, assume that $\hat{m}$ is $C^2$ in a neighborhood of $SO(2)$ and  $\R^{2 \times 2}_{\rm sym} \subset {\rm ker}(D\hat{m}(\Id))$. Let ${\cal F}_\eps(u) = {\cal E}_\eps(u) + \frac{1}{\eps}\int_{\Omega_\eps} f_{\kappa}(\nabla \tilde{y})$ with 
\begin{align}\label{eq:fk-def}
f_{\kappa}(F) = \begin{cases}
\kappa(1- \e_1 \cdot \hat{m}(F)), & |F| \le T, \\ 0 & \text{else},
\end{cases}
\end{align}
for $F \in \R^{2 \times 2}$, where $T,\kappa>0$.  Likewise, we define ${\cal F}^\chi_\eps$. In Lemma \ref{lemma: W,f} below we show that $W_{\Delta,\chi}(F) + f_{\kappa}(F) \ge C|F - \Id|^2$ for all $F \in \R^{2 \times 2}$ with $|F| \le T$.

This implies that the sequence $({\cal F}^\chi_\eps)_\eps$ is equicoercive:  Given a sequence of displacement fields $(u_\eps)_\eps$ with ${\cal F}^\chi_\eps(u_\eps)  + \Vert u_\eps \Vert_\infty \le C$ we find a subsequence converging in the sense of Definition \ref{def: convergence}. Indeed, we get that $\# {\cal C}^*_\eps \le \frac{C}{\eps}$, where ${\cal C}^*_\eps := \lbrace \Delta \in \tilde{\cal C}_\eps: |(\Id + \sqrt{\eps} \tilde{u}_\eps)_\Delta| > T \rbrace$. By  Lemma \ref{lemma: W,f} we then get $\| \nabla \tilde{u}_{\eps} \|_{L^2(\tilde{\Omega} \setminus \cup_{\triangle \in {\cal C}^*_{\eps}} \triangle)} \le C$  and therefore condition (ii) in Definition \ref{def: convergence} is satisfied. By  an $SBV$ compactness theorem (see \cite{Ambrosio-Fusco-Pallara:2000}) we then find a (not relabeled) subsequence such that $\tilde{u}_\eps \chi_{\tilde{\Omega}_\eps  \setminus \cup_{\triangle \in {\cal C}^*_{\eps}} \triangle} \to u$ in $L^1$ for some $u\in SBV^2(\tilde{\Omega})$. This together with $\Vert u_\eps \Vert_\infty \le C$ and $|\bigcup_{\triangle \in {\cal C}^*_{\eps}} \triangle| \le C\eps$ implies that also condition (i) in Definition \ref{def: convergence}(i) holds with this function $u$. 

Define $\hat{m}_1: \R^{2 \times 2} \to [-1,1]$ by $\hat{m}_1 = \e_1 \cdot \hat{m}$ and let $\hat{Q} = D^2\hat{m}_1(\Id)$ be the Hessian at the identity. We introduce the limiting functional ${\cal F}: SBV^2(\tilde{\Omega}) \to [0,\infty)$ given by
\begin{equation*}
{\cal F}(u) =  {\cal E}(u) - \frac{\kappa}{2}\int_{\Omega} \hat{Q}(\nabla u).
\end{equation*}
We then obtain a $\Gamma$-convergence result similar to Theorem \ref{th: Gamma-convergence}.

\begin{theorem}\label{th: Gamma-convergence2}
The assertions of Theorem \ref{th: Gamma-convergence} remain true when ${\cal E}_\eps$, ${\cal E}^\chi_\eps$ and ${\cal E}$ are replaced by ${\cal F}_\eps$, ${\cal F}^\chi_\eps$ and ${\cal F}$, respectively.
\end{theorem}

\subsection{Analysis of a limiting cleavage problem}

We now analyze the limiting functional ${\cal E}$ for a rectangular slab $\Omega = (0, l) \times (0, 1)$ with $l \ge \frac{1}{\sqrt{3}}$ under uniaxial extension in $\e_1$ direction. We determine the minimizers and prove uniqueness up to translation of the specimen and the crack line for the boundary conditions 
\begin{align}\label{eq:BC}
  u_1 = 0 \text{ for } x_1 = 0 \quad \text{ and } \quad u_1= a l \text{ for } x_1 = l. 
\end{align} 
 (More precisely: $u \in SBV^2((-\eta, l+ \eta) \times (0, 1))$ with $u_1(x) = 0$ for $x \le 0$ and $u_1(x) = al$ for $x \ge l$.) 
Note that we can investigate the limiting problem without any assumption on the second component of the boundary displacement. Let $\gamma = \max\{|\vv_1 \cdot \e_2|, |\vv_2 \cdot \e_2|, |(\vv_2 - \vv_1) \cdot \e_2|\}$ and $\vv_{\gamma} \in {\cal V}$ such that $\gamma = |\vv_{\gamma} \cdot \e_2|$. We note that $\gamma$ takes values in $[\frac{\sqrt{3}}{2},1]$ and $\vv_{\gamma}$ is unique if and only if $\phi \ne 0$. It turns out that the specimen shows perfect elastic behavior up to the critical boundary displacement 
$$ a_{\rm crit} = \sqrt{\frac{2\sqrt{3}\beta}{\alpha \gamma l}}. $$ 
Beyond critical loading the body fails by breaking into two pieces. 

\begin{theorem}\label{th: sharpness, uniqueness}
Let $a \neq a_{\rm crit}$. Then 
$$ \min \big\{ {\cal E}(u) : u \text{ satisfies \eqref{eq:BC}} \big\} 
   = \min \left\{\frac{ \alpha l}{\sqrt{3}} a^2 , \frac{2 \beta}{\gamma} \right\}. $$ 
All minimizers of ${\cal E}$ subject to \eqref{eq:BC} are of the following form:  
\begin{itemize} 
\item[(i)] If $a < a_{\rm crit}$, then 
$$ u^{\rm el}(x)  
   = (0,s) + \begin{pmatrix}
a & 0 \\ 0 & -\frac{a}{3}
\end{pmatrix} x  $$ 
for some $s \in \R$. 
\item[(ii)] If $a > a_{\rm crit}$ and $\phi \neq 0$ then 
$$ u^{\rm cr}(x)
  = \begin{cases} 
      (0,s) & \mbox{\rm for } x \mbox{\rm  \ to the left of } (p,0) + \R\vv_{\gamma}, \notag \\ 
      (a l, t) & \mbox{\rm for } x \mbox{\rm  \ to the right of } (p,0) + \R\vv_{\gamma}, 
    \end{cases} $$ 
for some $s,t \in \R$ and $p \in (0,l)$ such that $(p,0) + \R\vv_{\gamma}$ intersects both the segments $(0, l) \times \{0\}$ and $(0, l) \times \{1\}$. 
\item[(iii)] If $a > a_{\rm crit}$ and $\phi = 0$ then 
$$ u^{\rm cr}(x)
  = \begin{cases} 
      (0,s) & \mbox{\rm if } 0 < x_1 < h(x_2), \notag \\ 
      (a l, t) & \mbox{\rm if } h(x_2) < x_1 < l, 
    \end{cases} $$ 
for a Lipschitz function $h:(0,1) \to [0,l]$ with $|h'| \le \frac{1}{\sqrt{3}}$ a.e. and constants $s,t \in \R$. 
\end{itemize}
\end{theorem}

This theorem will be addressed in Section \ref{sec:convergence-variational-problem}. An analogous result for isotropic, incompressible materials has been obtained recently by Mora-Corral \cite{Mora:2010}. Theorem \ref{th: sharpness, uniqueness} is an extension of this result to anisotropic, compressible brittle materials in the framework of linearized elasticity. 

In particular, as mentioned above we see that all the optimal configurations show purely elastic behavior in the subcritical case and complete fracture in the supercritical regime. The crack minimizer in (ii) for $\phi \ne 0$ is broken parallel to $\R \vv_{\gamma}$ which proves that cleavage occurs along crystallographic lines, while in the symmetric case $\phi = 0$ cleavage in general fails.

\subsection{Applications: Cleaved crystals and fractured magnets}

As applications of the converging results for the energy functionals ${\cal E}^\chi_\eps$ and ${\cal F}^\chi_\eps$ we consider cleaved crystals and fractured magnets, respectively. In the first model a mild equicoercivity of the sequence $({\cal E}^\chi_\eps)_\eps$ is guaranteed by investigating a specific boundary value problem, in the latter model an external field provides an even stronger equicoercivity condition.

\subsubsection{Uniaxially strained crystals}

Theorem \ref{th: Gamma-convergence} in combination with Theorem \ref{th: sharpness, uniqueness} gives a new perspective to the results on crystal cleavage of \cite{FriedrichSchmidt:2014}. Let $\Omega = (0, l) \times (0, 1)$ with $l \ge \frac{1}{\sqrt{3}}$. For $\tilde{\Omega} = (-\eta, l + \eta) \times (0,1)$  and  $a \ge 0$ set  
\begin{align*} 
\begin{split}
  {\cal A}( a) 
  = \big\{ u = &(u_1, u_2) : \eps {\cal L} \cap \tilde{\Omega} \to \R^2 :  \\  & u(x) = g(x) \text{ for } x_1 \le \eps \text{ and } 
   x_1 \ge l - \eps 
     \text{ for  some } g \in {\cal G}(a)\big\},  
     \end{split}
\end{align*}
where ${\cal G}(a) := \lbrace g \in W^{1,\infty}(\tilde{\Omega}): g_1(x) =  0 \text{ for } x_1 \le \eps, ~ g_1(x) =  al  \text{ for } x_1 \ge l - \eps \rbrace$. In \cite[Theorem 2.1]{FriedrichSchmidt:2014} we proved that the limiting minimal energy leads to a universal cleavage law of the form 
\begin{equation*}
  \lim_{\eps \to 0} \inf \left\{ {\cal E}_{\eps}( u) : u \in {\cal A}(a) \right\} 
  = \min \left\{\frac{ \alpha l}{\sqrt{3}} a^2, \frac{2 \beta}{\gamma} \right\},  
\end{equation*}
independent of the particular shape of the interatomic potential $W$. Optimal configurations are given by the constant sequences $u_\eps = u^{\rm el}$ in the subcritical case $a \le a_{\rm crit}$ and $u_\eps = u^{\rm cr}$ in the supercritical case $a \ge a_{\rm crit}$, respectively, with $u^{\rm el}$ and $u^{\rm cr}$ as in Theorem \ref{th: sharpness, uniqueness}.

In fact, the above given configurations provide a characterization of all minimizing sequences in the sense that, all low energy sequences $(u_\eps)_\eps$ satisfying
\begin{align}\label{eq:almost-min}
{\cal E}^\chi_\eps(u_\eps) = \inf \lbrace {\cal E}^\chi_\eps(u): u \in {\cal A}(a) \rbrace + O(\eps)
\end{align}
and $\sup_{\eps} \|  u_{\eps} \|_{\infty}< \infty$ converge--up to subsequences--in the sense of Definition \ref{def: convergence} to $u^{\rm el}$ if $a < a_{\rm crit}$ or $u^{\rm cr}$ if $a > a_{\rm crit}$ for suitable $s, t, p$ and $g$, respectively. This is a direct consequence of \cite[Theorem 2.3 and Corollary 2.4]{FriedrichSchmidt:2014}. (The convergence obtained in \cite{FriedrichSchmidt:2014} is even stronger.)

One implication of \cite[Theorem 2.3 and Corollary 2.4]{FriedrichSchmidt:2014} is that, under the tensile boundary conditions $ u_{\eps} \in {\cal A}(a)$, the requirement that $ u_{\eps}$ be an almost energy minimizer satisfying \eqref{eq:almost-min}, guarantees the existence of a  subsequence converging in the sense of Definition \ref{def: convergence}. In particular, the sequence $({\cal E}^\chi_{\eps})$ is mildly equicoercive. A fundamental theorem of $\Gamma$-convergence (see, e.g., \cite[Theorem 1.21]{Braides:02}) implies that such low energy sequences converge to limiting configurations $u^{\rm el}$, respectively, $u^{\rm cr}$, in the sense of Definition \ref{def: convergence}. Consequently, in this way we have re-derived the convergence result \cite[Corollary 2.4]{FriedrichSchmidt:2014} (in the sense of Definition \ref{def: convergence}).

\subsubsection{Permanent magnets in an external field}
 
Assume that the material is a permanent magnet and let $\e_1$ be the magnetization direction. We suppose that there is a constitutive relation between $\nabla \tilde{y}(x)$ and the local magnetization direction $\hat{m}(\tilde{y},x) \in S^1$ of the deformed configuration $\tilde{y}$ at some point $x \in \Omega$, which is of the form $\hat{m}(\tilde{y},x) = \hat{m} (\nabla \tilde{y}(x))$ with $\hat{m}$ as defined in Section \ref{sec: con}. Let $H_{\rm ext}: \R^2 \to \R^2$ be an external magnetic field. The magnetic energy corresponding to the deformation $y = \id + \sqrt{\eps} u$ is then given by
$${\cal E}^{\rm mag}_\eps(u) = -\frac{1}{\eps} \int_{\Omega_\eps} H_{\rm ext} \cdot \hat{m}(\nabla \tilde{y}),$$
i.e. alignment of the magnetization direction with the external field is energetically favored. The total energy of the system is given by
$${\cal E}^{\rm tot}_\eps = {\cal E}^\chi_\eps + {\cal E}^{\rm mag}_\eps.$$ 
We now suppose that the external field is homogeneous and satisfies without restriction $H_{\rm ext} = \kappa \e_1$ for $\kappa >0$. We then see that 
$${\cal F}_\eps = {\cal E}^{\rm tot}_\eps -  \frac{\kappa}{\eps}|\Omega_\eps| $$
with $f_{\kappa}$ as in \eqref{eq:fk-def} and corresponding ${\cal F}_\eps$. By Theorem \ref{th: Gamma-convergence2} we get  that ${\cal E}^{\rm tot}_\eps$ $\Gamma$-converges to ${\cal E}^{\rm tot} = {\cal F}$ after renormalization with the sequence of constants $(\frac{\kappa}{\eps}|\Omega_\eps|)_{\eps}$. (Obviously, a configuration minimizes ${\cal E}^{\rm tot}_\eps$ if and only if it minimizes ${\cal F}_{\eps}$.)

We consider a boundary value problem $\min_{u \in {\cal A}_g} {\cal E}^{\rm tot}(u)$ for $g \in W^{1,\infty}(\tilde{\Omega})$. Since the sequence $({\cal E}^{\rm tot}_\eps)_\eps$ is equicoercive as discussed in Section \ref{sec: con}, the theory of $\Gamma$-convergence implies $\lim_{\eps \to 0 }\min_{u \in {\cal A}_g} {\cal F}_\eps (u) = \min_{u \in {\cal A}_g} {\cal F}^{\rm tot}(u)$ and also convergence of the corresponding (almost) minimizers  of ${\cal E}^{\rm tot}_{\eps}$ in the sense of Definition \ref{def: convergence} is guaranteed. In this context, note that by a truncation argument taking $g\in W^{1,\infty}(\tilde{\Omega})$ into account, we may indeed assume that a low energy sequence satisfies $\sup_\eps \Vert u_\eps \Vert_\eps < + \infty$.

\section{Convergence of the variational problems}\label{sec:convergence-variational-problem}

\subsection{Preparations}\label{subsec:auxiliary-lemmas}

The goal of this section is the derivation of the $\Gamma$-convergence result for ${\cal E}_{\eps}$. We first collect some properties of the cell energy $W_{\triangle}$ proven in \cite[Section 3]{FriedrichSchmidt:2014} provided that $W$ satisfies the assumptions (i), (ii) and (iii).

\begin{lemma}\label{lemma:W-triangle-properties}  
$W_{\triangle}$ is 
\begin{itemize}
\item[(i)] frame indifferent: $W_{\triangle}(QF) = W_{\triangle}(F)$ for all $F \in \R^{2 \times 2}$, $Q \in O(2)$, 
\item[(ii)] non-negative and satisfies $W_{\triangle}(F) = 0$ if and only if $F \in O(2)$ and 
\item[(iii)] $\liminf_{|F| \to \infty} W_{\triangle}(F) = \liminf_{|F| \to \infty} W_{\triangle,\chi}(F) = \beta$. 
\end{itemize}
\end{lemma}

\begin{lemma}\label{lemma:W-linearized-properties}
Let $F = \Id + G$ for $G \in \R^{2 \times 2}$. Then for $|G|$ small 
$$ W_{\triangle}(F) = \frac{1}{2} Q(G) + o(|G|^2),$$
where $Q(G) = \frac{3 \alpha}{16} \left( 3 g_{11}^2 + 3 g_{22}^2 + 2 g_{11} g_{22} + 4 \left(\frac{g_{12}+g_{21}}{2} \right)^2 \right)$. 

In particular, $Q(G)$ only depends on the symmetric part $\left(G^{T} + G\right)/2$ of $G$. $Q$ is positive semidefinite and thus convex on $\R^{2 \times 2}$ and positive definite and strictly convex on the subspace $\R^{2 \times 2}_{\rm sym}$ of symmetric matrices. 
\end{lemma}

The following lemma provides useful lower bounds for the energy $W_{\triangle}$ and the pair interaction potential $W$. 

\begin{lemma}\label{lemma:quadratic lower bound} 
For all $T>1$ one has:
\begin{itemize}
\item[(i)] There exists some $c > 0$ such that $c \dist^2(F,O(2)) \leq  W_{\triangle} (F)$ for all $F \in \R^{2 \times 2}$ satisfying $|F|\leq T$.
\item[(ii)] For $\rho > 0$ there is an increasing, subadditive function $\psi^{\rho}:[0,\infty) \to (0, \infty)$ which satisfies $\psi^{\rho}(r) - \rho \leq W(r+1)$ for all $r\geq 0$ and $\psi(r) = \beta$ for all $r \geq c_{\rho}$ for some constant $c_{\rho}$ only depending on $\rho$.
\end{itemize}
\end{lemma}

\Proof
(i) This essentially follows from the expansion given in Lemma \ref{lemma:W-linearized-properties}. For details we refer to \cite[Lemma 3.5]{FriedrichSchmidt:2014}. 

(ii) We define 
$$ \bar{\psi}(r) 
= \begin{cases} 
                \eta r
          & \text{for } 0 \le r \le \frac{\beta}{\eta} , \\ 
        \beta
          & \text{for } r \ge \frac{\beta}{\eta}, 
     \end{cases} $$ 
for some $\eta > 0$ (depending on $\rho$) such that $\bar{\psi} - \rho \leq W$. Then we set $\psi^{\rho}(r) = \bar{\psi}(r+1)$. As $\psi^{\rho}$ is a concave function with $\psi^{\rho}(0) > 0$, it is subadditive. 
\eop

Moreover, we provide a lower bound for $W_{\Delta,\chi}(F) + f_{\kappa}(F)$ which implies the equicoercivity of  $({\cal F}^\chi_\eps)_\eps$.

\begin{lemma}\label{lemma: W,f}
Let $T > \sqrt{2}$. Then there are constants $C_1,C_2>0$ such that for all $F \in \R^{2 \times 2}$ with $|F| \le T$ we obtain
\begin{itemize}
\item[i)] $|\hat{m}(F) - \hat{m}(R(F))| \le C_1|F - R(F)|^2$, where $R(F) \in SO(2)$ is a solution of $ |F - R(F)| = \min_{R \in SO(2)}|F - R|$,
\item[(ii)] $W_{\Delta,\chi}(F) + f_{\kappa}(F) \ge C_2|F - \Id|^2$.
\end{itemize}
\end{lemma}

\Proof (i) Without restriction we may assume that $|F - R(F)|$ is small as otherwise the assertion is clear. So in particular, $R(F)$ is uniquely determined. Moreover, it suffices to consider $F \in \R^{2 \times 2}_{\rm sym}$ and $R(F) = \Id$. Indeed, once this is proved, we find $|\hat{m}(F) - \hat{m}(R(F))| = |R(F)\hat{m}(R(F)^T F) - R(F)\hat{m}(\Id))| \le C|R(F)^T F - \Id|^2$, as desired. 

Let $F \in \R^{2 \times 2}_{\rm sym}$, $R(F) = \Id$ and set $G = F - \Id$  with $G \in \R^{2 \times 2}_{\rm sym}$ small. As $\hat{m}$ is $C^2$ in a neighborhood of $SO(2)$ we derive $|\hat{m}(F)- \hat{m}(\Id)| \le |D\hat{m}(\Id) \, G| + C|G|^2 = C|G|^2$ as $\R^{2 \times 2}_{\rm sym} \subset {\rm ker}(D\hat{m}(\Id))$.

(ii) By Lemma \ref{lemma:quadratic lower bound}(i) the assertion is clear for all $|F| \le T$ with $c_0 \le \dist(F,O(2))$ for $c_0 >0$ and $C=C(c_0,T)$ sufficiently small. Otherwise, we again apply Lemma \ref{lemma:quadratic lower bound}(i) to obtain for $c_0$ small enough
$$W_{\Delta,\chi}(F) \ge C\dist^2(F,O(2)) + \chi(F) \ge C\dist^2(F,SO(2)) = C|F - R(F)|^2.$$
For convenience we write $r_{ij} = \e^T_i R(F) \e_j$ for $i,j = 1,2$. As $r_{12}^2 = r_{21}^2 = 1- r_{11}^2$ we find $1- r_{11} = 1- r^2_{11} + r_{11}(r_{11}-1) = r_{12}^2 +  (1- r_{11})^2  -(1- r_{11})$. Thus, recalling $\hat{m}(R) = R\e_1$ for all $R \in SO(2)$ and applying (i) we get for $0 < c \le \kappa$ small enough 
\begin{align}\label{eq: r_11}
\begin{split}
W_{\Delta,\chi}&(F) + f_{\kappa}(F) \\
& \ge C|F - R(F)|^2 + c(1 - \e_1 \cdot \hat{m}(R(F))) + c\e_1\cdot(\hat{m}(R(F)) - \hat{m}(F)) \\ &
\ge C|F - R(F)|^2 + c(1 - \e_1^T R(F) \e_1) - c C_{1}|F - R(F)|^2 \\ 
& \ge  \frac{C}{2} |F - R(F)|^2 + \frac{c}{2} (1 - r_{11})^2 + \frac{c}{2} r_{12}^2 
\ge C_{2} |F - \Id|^2, 
\end{split}
\end{align}
as desired. \eop 

As a further preparation we modify the interpolation $\tilde{y}$ on triangles with large deformation: We fix a threshold explicitly as $R = 7$ and let $\bar{{\cal C}}_{\eps} \subset \tilde{{\cal C}}_{\eps}$ be the set of those triangles where $|(\tilde{y})_{\triangle}| > R$. By definition of the boundary values in \eqref{eq: Ag} we find $\bar{{\cal C}}_{\eps} \subset {{\cal C}}_{\eps}$ for $\eps$ small enough. We introduce another interpolation $y'$ which leaves $\tilde{y}$ unchanged on $\triangle \in \tilde{\cal C}_{\eps} \setminus \bar{{\cal C}}_{\eps}$ and replaces $\tilde{y}$ on $\triangle \in \bar{{\cal C}}_{\eps}$ by a discontinuous function with constant derivative satisfying $|(y')_{\triangle}| \leq R$. In fact, by introducing jumps we achieve a release of the elastic energy. Note that $y' \in SBV (\tilde{\Omega}_\eps)$. 

More precisely, note that on $\triangle \in \bar{{\cal C}}_{\eps}$ we have $|(\tilde{y})_{\triangle} \, \vv| \geq 2$ for at least two springs $\vv \in {\cal V}$. Indeed, using the elementary identity 
\begin{align*}
  \sum_{\vv \in {\cal V}} \langle \vv, H \vv \rangle^2 
  = \frac{3}{8} \left( 2\trace(H^2) + (\trace H)^2 \right) 
  \ge \frac{3}{8} (\trace H)^2 
\end{align*}
for any $H \in \R^{2 \times 2}_{\rm sym}$, we find that $|F| > 7$ implies 
\begin{align*}
  \sum_{\vv \in {\cal V}} |F \vv|^4 
  = \sum_{\vv \in {\cal V}} \langle \vv, F^T F \vv \rangle^{2} 
  \ge \frac{3}{8} (\trace (F^T F))^2 
  = \frac{3}{8} |F|^4 
\end{align*} 
and so $\max_{\vv \in {\cal V}} |F \vv|^{4} > \frac{7^4}{8} > 4^4$. Hence, $|F \vv| > 4$ for at least one $\vv \in {\cal V}$ and at least two springs are elongated by a factor larger than $2$. For $m=2,3$ let $\bar{{\cal C}}_{\eps,m} \subset \bar{{\cal C}}_{\eps}$ be the set of triangles where $|(\tilde{y})_{\triangle} \, \vv| \geq 2$ holds for exactly $m$ springs $\vv \in {\cal V}$. For $i,j,k=1,2,3$ pairwise distinct let $h_{i}$ denote the segment between the centers of the sides in $\vv_{j}$ and $\vv_{k}$ direction and define the set $V_{i}= h_{j} \cup h_{k}$.

We now construct $y' \in SBV^2(\tilde{\Omega}_{\eps})$. On $\triangle \in \tilde{\cal C}_{\eps} \setminus \bar{{\cal C}}_{\eps}$ we simply set $y' = \tilde{y}$. On $\triangle \in \bar{{\cal C}}_{\eps,2}$, assuming $|(\tilde{y})_{\triangle} \, \vv_i|\leq 2$, we choose $y'$ such that $\nabla y'$ assumes the constant value $(y')_{\triangle}$ on $\triangle$ with $(y')_{\triangle} \, \vv_i =  (\tilde{y})_{\triangle} \, \vv_i$ and $|(y')_{\triangle} \, \vv| = 1$ for $\vv \in {\cal V}\setminus \left\{\vv_i\right\}$. Moreover, we ask that $y' = \tilde{y}$ at the three vertices and on the side orientated in $\vv_i$ direction. This can and will be done in such a way that $y'$ is continuous on $\text{int}(\triangle)\setminus h_{i}$. We note that the definition of $(y')_{\triangle}$ is unique up to a reflection, unless $(\tilde{y})_{\triangle} \vv_{i} = 0$. We may and will assume that 
\begin{equation}\label{eq: SO, O}
\dist\left((y')_{\triangle}, SO(2)\right) \leq \dist\left((y')_{\triangle}, O(2) \setminus SO(2)\right).
\end{equation}
For $\triangle \in \bar{{\cal C}}_{\eps,3}$ we set $(y')_{\triangle} = \Id$ and $y' = \tilde{y}$ at the three vertices such that $y'$ is continuous on $\text{int}(\triangle)\setminus V_{i}$ for some $i \in \{1,2,3\}$. Here, the index $i$ can be taken arbitrarily at first. However, in what follows it will also be necessary to use the following unambiguously defined `variants' of $y'$: If on every $\triangle \in \bar{{\cal C}}_{\eps,3}$ the set $V_{i}$ is chosen as the jump set of $y'$ we denote this interplation explicitly as $y'_{V_{i}}$. 

We define the interpolation $u'$ for the rescaled displacement field by $u' = \frac{1}{\sqrt{\eps}} (y' - \id)$. We note that by construction also on an edge $[p,q] \subset \partial \triangle$ for $\triangle \in \bar{{\cal C}}_{\eps}$ jumps may occur. There, however, the jump height $|[u_{\eps}]|$ can be bounded by  
\begin{equation}\label{eq: jump height, rest jumps}
  |[u'_{\eps}](x)| 
  \leq \eps \left\|\nabla u'_{\eps}\right\|_{\infty} \leq \eps \cdot c\eps^{-\frac{1}{2}} = c\sqrt{\eps}
\end{equation}
for a constant $c > 0$ independent of $\eps$ and $x \in [p, q]$. This holds since the interpolations are continuous at the vertices.

The following lemma shows that we may pass from $\tilde{u}_{\eps}$ to $u'_{\eps}$ without changing the limit.
\begin{lemma}\label{lemma: tilde'}
If $u_{\eps} \to u$ in the sense of Definition \ref{def: convergence} and ${\cal E}(u_{\eps})$ is uniformly bounded, then $\chi_{\tilde{\Omega}_{\eps}} u'_{\eps} \to u$ in $L^{1}(\tilde{\Omega})$, $\chi_{\tilde{\Omega}_{\eps}}\nabla u_{\eps}' \weakly \nabla u$ in $L^2(\tilde{\Omega})$ and ${\cal H}^1(J_{u'_{\eps}})$ is uniformly bounded.   
\end{lemma}

\Proof We first note that there is some $M > 0$ such that
\begin{equation}\label{eq: jump set bound}
\#\bar{{\cal C}}_{\eps} \leq \frac{M}{\eps} 
\end{equation}
for all $\eps >0$. To see this, we just recall that every triangle $\triangle \in \bar{{\cal C}}_{\eps}$ provides at least the energy $\eps \inf\left\{W(r): r \geq 2\right\}$. In fact we may assume that ${\cal C}^*_{\eps} = \bar{\cal C}_{\eps}$ in Definition  \ref{def: convergence} as for $\Delta \in {\cal C}^*_{\eps} \setminus \bar{\cal C}_{\eps}$ we have $|(\tilde{u}_{\eps})_{\triangle}| \le \frac{C}{\sqrt{\eps}} |(\tilde{y}_{\eps})_{\triangle} - \Id| \le \frac{C}{\sqrt{\eps}}$ and so 
\begin{align*} 
  \| \nabla \tilde{u}_{\eps} \|_{L^2(\tilde{\Omega}_{\eps} \setminus \cup_{\triangle \in \bar{\cal C}_{\eps}} \triangle)} 
  &\le \| \nabla \tilde{u}_{\eps} \|_{L^2(\tilde{\Omega}_{\eps} \setminus \cup_{\triangle \in {\cal C}^*_{\eps}} \triangle)} 
    + \| \nabla \tilde{u}_{\eps} \|_{L^2(\cup_{\triangle \in {\cal C}^*_{\eps} \setminus \bar{\cal C}_{\eps}} \triangle)} \\ 
  & \le C + \left( \# ( {\cal C}^*_{\eps} \setminus \bar{\cal C}_{\eps} ) \frac{\sqrt{3}\eps^2}{4} \cdot \frac{C}{\eps} \right)^{\frac{1}{2}} 
   \le C. 
\end{align*}
It follows that $\chi_{\tilde{\Omega}_{\eps}}\nabla u_{\eps}'$ is bounded uniformly in $L^2$ and, in particular, equiintegrable. Finally, the jump lengths ${\cal H}^1(J_{u'_{\eps}})$ are readlily seen to be bounded by $C \eps \#\bar{{\cal C}}_{\eps} \le C$. But then Ambrosio's compactness Theorem for GSBV \cite[Theorem 2.2]{Ambrosio:90} shows that indeed $\chi_{\tilde{\Omega}_{\eps}}\nabla u_{\eps}' \weakly \nabla u$ in $L^2(\tilde{\Omega})$. \eop

\subsection{The $\Gamma$-$\liminf$-inequality}

With the above preparations at hand, we may now prove the $\Gamma$-$\liminf$-inequality in Theorem \ref{th: Gamma-convergence}. \smallskip

\noindent {\em Proof of Theorem \ref{th: Gamma-convergence}(i).} 
Let $(g_\eps)_\eps \in W^{1,\infty}(\tilde{\Omega})$ with $ \sup_{\eps} \Vert g_\eps\Vert_{W^{1,\infty}(\tilde{\Omega})} < + \infty$ be given. Let $u \in SBV^2(\tilde{\Omega})$ and consider a sequence $u_{\eps} \subset SBV^2(\tilde{\Omega}_{\eps})$ with $u_{\eps} \in {\cal A}_{g_{\eps}}$ converging to $u$ in $SBV^2$ in the sense of Definition \ref{def: convergence}. We split up the energy into bulk and crack parts neglecting the contribution $\eps E^{\rm boundary}_{\eps}$ from the boundary layers:
\begin{align}\label{eq: auxiliary lower bound}
\begin{split}
{\cal E}_{\eps}(u_{\eps}) 
   & \geq \eps \sum_{\triangle  \in {\cal C}_{\eps} \setminus \bar{{\cal C}}_{\eps}} W_{\triangle}((\tilde{y}_{\eps})_{\triangle}) + \eps \sum_{\triangle \in \bar{{\cal C}}_{\eps}} W_{\triangle}((\tilde{y}_{\eps})_{\triangle})\\
   & = \frac{4}{\sqrt{3}\eps} \int_{\Omega_{\eps}} W_{\triangle}\left(\Id + \sqrt{\eps} \nabla u'_{\eps}\right) + \eps \sum_{\triangle \in \bar{{\cal C}}_{\eps}} \ \sum_{\vv \in {\cal V}, \atop |(\tilde{y}_{\eps})_{\triangle} \, \vv|>2} \frac{1}{2} W\left(|(\tilde{y}_{\eps})_{\triangle} \, \vv|\right)\\
   & =: {\cal E}^{\rm elastic}_{\eps} (u_{\eps}) + {\cal E}^{\rm crack}_{\eps} (u_{\eps}).
   \end{split}
\end{align}
We note that by contruction of the interpolation $u'_{\eps}$ we may take the integral over $\Omega_{\eps}$. As both parts separate completely in the limit, we discuss them individually. \smallskip

\noindent\textit{Elastic energy.} 
We first concern ourselves with the elastic part of the energy. We recall $W_{\triangle}(\Id + G) = \frac{1}{2}Q(G) + \omega(G)$ with $\sup \left\{\frac{\omega(F)}{|F|^2} : |F|\leq \rho\right\} \rightarrow 0$ as $\rho \rightarrow 0$. Let $\chi_{\eps}(x) := \chi_{[0,\eps^{-1/4})}(|\nabla u'_{\eps}(x)|)$. Note that for $F \in \R^{2 \times 2}$, $r>0$ one has $Q(rF) = r^2 Q(F)$. We compute 
\begin{align*}
  {\cal E}^{\rm elastic}_{\eps} (u_{\eps}) 
  & \geq \frac{4}{\sqrt{3}} \int_{\Omega_{\eps}} \chi_{\eps}(x)\left( \frac{1}{2} Q( \nabla u'_{\eps}) + \frac{1}{\eps} \omega\left(\sqrt{\eps}\nabla u'_{\eps}(x)\right) \right) \, dx. 
\end{align*}
The second term of the integral can be bounded by 
$$\chi_{\eps} |\nabla u'_{\eps}|^2\frac{\omega\left(\sqrt{\eps} \nabla u'_{\eps} \right)}{|\sqrt{\eps}\nabla u'_{\eps} |^2}.$$
Since $\nabla u'_{\eps}$ is bounded in $L^2$ and $\chi_{\eps} \frac{\omega\left(\sqrt{\eps} \nabla u'_{\eps} \right)}{|\sqrt{\eps} \nabla u'_{\eps} |^2}$ converges uniformly to $0$ as $\eps \rightarrow 0$ it follows that 
\begin{align*}
\liminf_{\eps \rightarrow 0}{\cal E}^{\rm elastic}_{\eps} (u_{\eps}) & \geq \liminf_{\eps \rightarrow 0} \frac{4}{\sqrt{3}} \int_{\Omega_{\eps}} \chi_{\eps}(x) \frac{1}{2} Q( \nabla u'_{\eps}(x)) \, dx\\
& \geq \liminf_{\eps \rightarrow 0} \frac{4}{\sqrt{3}} \int_{\Omega}  \frac{1}{2} Q(\chi_{\Omega_{\eps}} \chi_{\eps}(x) \nabla u'_{\eps}(x)) \, dx.
\end{align*}
By assumption $\chi_{\Omega_{\eps}} \nabla u'_{\eps} \rightharpoonup \nabla u$ weakly in $L^2$. As $\chi_{\eps} \rightarrow 1$ boundedly in measure on $\Omega$, it follows $\chi_{\Omega_{\eps}}\chi_{\eps} \nabla u'_{\eps} \rightharpoonup u$ weakly in $L^2(\Omega)$. By lower semicontinuity (Q is convex by Lemma \ref{lemma:W-linearized-properties}) we conclude recalling that $Q$ only depends on the symmetric part of the gradient:
\begin{align*}
\liminf_{\eps \rightarrow 0}{\cal E}^{\rm elastic}_{\eps} (u_{\eps}) 
 \geq   \frac{4}{\sqrt{3}} \int_{\Omega}  \frac{1}{2} Q(e(u(x))) \, dx.
\end{align*}  
\smallskip

\noindent\textit{Crack energy.}
By construction the functions $u'_{\eps}$ have jumps on destroyed triangles $\triangle \in \bar{{\cal C}}_{\eps}$. We now write the energy of such a triangle in terms of the jump height $\left[u\right] = u^{+} - u^{-}$. We first concern ourselves with a triangle $\triangle \in \bar{{\cal C}}_{\eps,3}$. For the variant $u'_{\eps,V_{i}}$, $i=1,2,3$ we consider the springs in $\vv_{j},\vv_{k}$ direction for $j,k \neq i$. Thus, we compute 
\begin{equation}\label{eq: jump height}
\eps (\tilde{y}_{\eps})_{\triangle} \, \vv_{j} = \eps (y'_{\eps})_{\triangle} \, \vv_{j} + [y'_{\eps, V_{i}}]_{h_{\vv_{k}}} = \eps \vv_{j} + \sqrt{\eps} [u'_{\eps, V_{i}}]_{h_{\vv_{k}}},
\end{equation}
where $[u'_{\eps, V_{i}}]_{h_{\vv_{k}}}$ denotes the jump height on the set $h_{\vv_{k}}$. Here and in the following equations, the same holds true if we interchange the roles of $j$ and $k$. We claim that 
\begin{equation}\label{eq: jump height estimate}
  |(\tilde{y}_{\eps})_{\triangle} \, \vv_{j}| 
  \geq \eps^{\frac{1}{4}} \left|\frac{1}{\sqrt{\eps}} [u'_{\eps, V_{i}}]_{h_{\vv_{k}}}\right| + 1.
\end{equation}
Indeed, for $|\frac{1}{\sqrt{\eps}} [u'_{\eps, V_{i}}]_{h_{\vv_{k}}}| \leq \eps^{-\frac{1}{4}}$ this is clear since $|(\tilde{y}_{\eps})_{\triangle} \, \vv_{j}| \geq 2$. Otherwise, applying (\ref{eq: jump height}) we compute for $\eps$ small enough:
\begin{align*}
  |(\tilde{y}_{\eps})_{\triangle} \, \vv_{j}| 
  &= \left|\frac{1}{\sqrt{\eps}} [u'_{\eps, V_{i}}]_{h_{\vv_{k}}} + \vv_{j}\right| 
   \geq \left|\frac{1}{\sqrt{\eps}} [u'_{\eps, V_{i}}]_{h_{\vv_{k}}}\right| - 1  \\ 
  &\geq\eps^{\frac{1}{4}}\left|\frac{1}{\sqrt{\eps}} [u'_{\eps, V_{i}}]_{h_{\vv_{k}}}\right| + \left(1 - \eps^{\frac{1}{4}}\right)\eps^{-\frac{1}{4}} - 1 \\
  &= \eps^{\frac{1}{4}}\left|\frac{1}{\sqrt{\eps}} [u'_{\eps, V_{i}}]_{h_{\vv_{k}}}\right| - 2 + \eps^{-\frac{1}{4}} 
   \geq \eps^{\frac{1}{4}}\left|\frac{1}{\sqrt{\eps}} [u'_{\eps, V_{i}}]_{h_{\vv_{k}}}\right| + 1.
\end{align*} 
Let $\rho > 0$ sufficiently small. Applying Lemma \ref{lemma:quadratic lower bound}(ii) there is an increasing subadditive function $\psi^{\rho} 0$ with $\psi^{\rho}(r-1) - \rho \leq W(r)$ for $r \geq 1$. We define $\tilde{\psi}^{\rho} = \psi^{\rho} - \rho$. The monotonicity of $\psi^{\rho}$ and (\ref{eq: jump height estimate}) yield 
\begin{equation}\label{eq: jump energy estimate}
  W(|(\tilde{y}_{\eps})_{\triangle} \, \vv_{j}|) 
  \geq \tilde{\psi}^{\rho}(|(\tilde{y}_{\eps})_{\triangle} \, \vv_{j}|-1) 
  \geq \tilde{\psi}^{\rho}\left(\left|\eps^{-\frac{1}{4}} [u'_{\eps, V_{i}}]_{h_{\vv_{k}}}\right|\right).
\end{equation}
Now for $\triangle \in \bar{{\cal C}}_{\eps,3}$ we may estimate the energy as follows:
\begin{align*}
  W_{\triangle} \left( (\tilde{y}_{\eps})_{\triangle}\right) 
  &= \frac{1}{2} \sum^{3}_{l=1} W(|(\tilde{y}_{\eps})_{\triangle} \, \vv_{l}|) \\
  &\geq \frac{1}{4} \sum^{3}_{i=1} \left\{\tilde{\psi}^{\rho}\left(\eps^{-\frac{1}{4}}|[u'_{\eps, V_{i}}]_{h_{\vv_{k}}}|\right) + \tilde{\psi}^{\rho}\left(\eps^{-\frac{1}{4}}|[u'_{\eps, V_{i}}]_{h_{\vv_{j}}}|\right) \right\} =: W_{\triangle,3}\left( (\tilde{y}_{\eps})_{\triangle}\right), 
\end{align*}
where $i,j,k =1,2,3$ are pairwise distinct. With $\nu^{(i)}_u = \nu_{u'_{\eps, V_{i}}}$ we can also write
$$ W_{\triangle,3}\left( (\tilde{y}_{\eps})_{\triangle}\right) 
   = \frac{1}{4} \cdot \frac{2}{\eps} \cdot \frac{2}{\sqrt{3}} \sum^{3}_{i=1} \int_{h_{\vv_{j}} \cup h_{\vv_{k}}} \tilde{\psi}^{\rho}\left(\eps^{-\frac{1}{4}} |[u'_{\eps, V_{i}}]|\right) \left(|\vv_{j}\cdot \nu^{(i)}_{u}| + |\vv_{k}\cdot \nu^{(i)}_{u}|\right) d{\cal H}^1. $$
The factors in front occur since ${\cal H}^1(h_{\vv_{j}}) = \frac{\eps}{2}$ and, letting $\nu_{j}$ be a normal of $h_{\vv_{j}}$, one has $|\nu_{j} \cdot \vv_{j}| = 0$ and $|\nu_{j} \cdot \vv_{k}| = \frac{\sqrt{3}}{2}$. Consequently, defining $\phi^{\rho}_{i}(r,\nu) = \psi^{\rho}(r) \left(|\vv_{j}\cdot \nu| + |\vv_{k}\cdot \nu|\right)$ and $\tilde{\phi}^{\rho}_{i}(r,\nu) = \tilde{\psi}^{\rho}(r) \left(|\vv_{j}\cdot \nu| + |\vv_{k}\cdot \nu|\right)$, respectively, we get 
\begin{equation*}
  W_{\triangle,3}\left( (\tilde{y}_{\eps})_{\triangle}\right) 
  = \frac{1}{\sqrt{3}{\eps}} \sum^{3}_{i=1} \int_{J_{u'_{\eps, V_{i}}} \cap {\text{int}}(\triangle)}  \tilde{\phi}^{\rho}_{i}(\eps^{-\frac{1}{4}}|[u'_{\eps, V_{i}}]|,\nu^{(i)}_{u}) \, d{\cal H}^1
\end{equation*}
on every $\triangle \in \bar{{\cal C}}_{\eps,3}$. For $\triangle \in \bar{{\cal C}}_{\eps,2}$ we proceed analogously. Assuming $|(\tilde{y}_{\eps})_{\triangle} \, \vv_{i}| \leq 2$ we compute for the springs in $\vv_{j}, \vv_{k}$ direction (abbreviated by $\vv_{j,k}$) as in (\ref{eq: jump height})
\begin{equation}\label{eq: [2] triangle}
  \eps (\tilde{y}_{\eps})_{\triangle} \, \vv_{j,k} = \eps (y'_{\eps})_{\triangle} \, \vv_{j,k} + \sqrt{\eps}  [u'_{\eps}]_{h_{\vv_{i}}}. 
\end{equation}  
Note that in this case we do not have to take a special variant of $u'_{\eps}$ into account. Repeating the steps (\ref{eq: jump height estimate}) and (\ref{eq: jump energy estimate}) we find
\begin{align*}
\frac{1}{2} \left(W(|(\tilde{y}_{\eps})_{\triangle} \, \vv_{j}|) + W(|(\tilde{y}_{\eps})_{\triangle} \, \vv_{k}|) \right) 
\geq \tilde{\psi}^{\rho}\left(\eps^{-\frac{1}{4}}|[u'_{\eps}]_{h_{\vv_{i}}}|\right) 
=: W_{\triangle,2}\left( (\tilde{y}_{\eps})_{\triangle}\right).  
\end{align*}
Noting that $|\vv_j \cdot \nu_i| = |\vv_k \cdot \nu_i| = \frac{\sqrt{3}}{2}$, $|\vv_i \cdot \nu_i| = 0$ and that every of these terms occurs twice in the sum of the right hand side of the following formula, it is not hard to see that this energy satisfies the same integral representation formula as $W_{\triangle,3}$: 
\begin{equation*}
  W_{\triangle,2}\left( (\tilde{y}_{\eps})_{\triangle}\right) 
  = \frac{1}{\sqrt{3}{\eps}} \sum^{3}_{i=1} \int_{J_{u'_{\eps, V_{i}}} \cap {\text{int}}(\triangle)}  \tilde{\phi}^{\rho}_{i}(\eps^{-\frac{1}{4}}|[u'_{\eps, V_{i}}]|,\nu^{(i)}_{u}) \, d{\cal H}^1.
\end{equation*}
(Recall that the interpolation variant $u'_{\eps, V_{i}}$ and its crack normal $\nu^{(i)}_{u}$ do not depend on $i$ on $\triangle \in \bar{{\cal C}}_{\eps,2}$.) Let $\sigma > 0$. Note that $\bar{{\cal C}}_{\eps} \subset {{\cal C}}_{\eps}$ for $\eps$ sufficiently small as $\sup_{\eps} \Vert g_\eps\Vert_{W^{1,\infty}(\tilde{\Omega})} < + \infty$. Thus, the crack energy can be estimated by
\begin{align*}
{\cal E}^{\rm crack}_{\eps} (u_{\eps}) 
& \geq \frac{1}{\sqrt{3}}  \sum_{i} \int_{J_{u'_{\eps, V_{i}}} \cap \tilde{\Omega}_{\eps}}  \tilde{\phi}^{\rho}_{i}(\eps^{-\frac{1}{4}}|[u'_{\eps, V_{i}}]|,\nu^{(i)}_{u}) \, d{\cal H}^1 - E^{\rho}_{\eps, \cup\partial \triangle}\left(\tilde{y}_{\eps}\right) \\
& \geq \frac{1}{\sqrt{3}} \sum_{i} \int_{J_{u'_{\eps, V_{i}}} \cap \tilde{\Omega}_{\eps}} \left(\phi^{\rho}_{i}(\sigma^{-1}|[u'_{\eps, V_{i}}]|,\nu^{(i)}_{u}) - 2 \rho\right) \, d{\cal H}^1 - E^{\rho}_{\eps, \cup\partial \triangle}\left(\tilde{y}_{\eps}\right),
\end{align*}
where {$E^{\rho}_{\eps, \cup \partial \triangle}\left(\tilde{y}_{\eps}\right)$} compensates for the extra contribution provided by jumps lying on the boundary of some $\triangle \in \bar{{\cal C}}_{\eps}$. We will show that this term vanishes in the limit. 

Now by construction the $\phi^{\rho}_{i}(r, \nu)$, $i = 1,2,3$, are products of a positive, increasing and concave function in $r$ and a norm in $\nu$. Moreover, $u'_{\eps}$ and its variants converge to $u$ in $L^1$ with $\nabla u'_{\eps}$ bounded in $L^2$ and thus equiintegrable. By Ambrosio's lower semicontinuity Theorem \cite[Theorem 3.7]{Ambrosio:90} we obtain 
$$ \liminf_{\eps \to 0} {\cal E}^{\rm crack}_{\eps} (u_{\eps})  \geq 
\frac{1}{\sqrt{3}} \int_{J_{u}} \sum_{i} \phi^{\rho}_{i}(\sigma^{-1}|[u]|,\nu_{u})\,  d{\cal H}^1 - C M \rho- \limsup_{\eps \to 0} E^{\rho}_{\eps, \cup \partial \triangle}\left(\tilde{y}_{\eps}\right),$$
where we used that $\sup_{\eps} {\cal H}^1(J_{u'_{\eps}}) \leq CM$ for a constant $C > 0$ by (\ref{eq: jump set bound}). We recall that $\psi^{\rho}(r) \rightarrow \beta$ for $r \rightarrow \infty$. In the limit $\sigma \rightarrow 0$ this yields
\begin{equation}\label{eq: rho crack energy}
\liminf_{\eps \to 0} {\cal E}^{\rm crack}_{\eps} (u_{\eps})  \geq 
\frac{1}{\sqrt{3}} \int_{J_{u}} 2\beta \sum_{\vv \in {\cal V}} |\vv \cdot \nu_{u}|\, d{\cal H}^1 - CM\rho - \limsup_{\eps \to 0} E^{\rho}_{\eps, \cup \partial \triangle}\left(\tilde{y}_{\eps}\right).
\end{equation}
Taking (\ref{eq: jump height, rest jumps}) and (\ref{eq: jump set bound}) into account we compute
\begin{align*}
  \limsup_{\eps \to 0} \sum_{\triangle \in \bar{{\cal C}}_{\eps}} \int_{\partial \triangle}|
  \tilde{\psi}^{\rho}\left(\eps^{-\frac{1}{4}} |[u'_{\eps}]|\right)| 
  &\leq \lim_{\eps \to 0} CM \sup \left\{|\psi^{\rho}\left(r\right) - \rho|: \, r\leq \eps^{-\frac{1}{4}} \cdot c\eps^{\frac{1}{2}}\right\} \\ 
  &= CM \rho. 
\end{align*}
This proves $\limsup_{\eps} |E^{\rho}_{\eps, \cup \partial \triangle}\left(\tilde{y}_{\eps}\right)| \leq \tilde{C}M \rho$ for some $\tilde{C} > 0$. We finally let $\rho \rightarrow 0$ in \eqref{eq: rho crack energy}. This finishes the proof of (i). \eop

We now prove the $\Gamma$-$\liminf$-inequality in Theorem \ref{th: Gamma-convergence2}. \smallskip
 
\noindent {\em Proof of Theorem \ref{th: Gamma-convergence2}, first part.} 
Following the proof of Theorem \ref{th: Gamma-convergence}(i) it suffices to show 
$$\liminf_{\eps \to 0}  \frac{1}{\eps}\int_{\Omega_\eps}\chi_\eps f_{\kappa}(\nabla y'_\eps) \ge - \frac{\kappa}{2}\int_\Omega \hat{Q}(\nabla u),$$
where $\hat{Q} = D^2 \hat{m}_1(\Id)$. Let $u'_\eps = \frac{1}{\sqrt{\eps}} (y'_\eps - \id)$. With a slight abuse of notation we set $e(F) = \frac{1}{2}(F^T + F)$ and $a(F) = F - e(F)$ for matrices $F \in \R^{2 \times 2}$. Let $F = \Id + \sqrt{\eps}G$ for $G \in \R^{2 \times 2}$. Linearization around the identity matrix yields $\dist(F,SO(2)) =  \sqrt{\eps}|e(G)| + \eps O(|G|^2)$. It is not hard to see that this implies 
\begin{align}\label{eq: Linearization}
R(F) = \Id + \sqrt{\eps} a(G) + \eps O(|G|^2),
\end{align}
where $R(F) \in SO(2)$ is as defined in Lemma \ref{lemma: W,f}.  As $\hat{m}(\Id) =\e_1$ and $e(G) \in {\rm ker}(D\hat{m}(\Id))$, we find by expanding $\hat{m}_1$
\begin{align}\label{eq: m-lin} 
\hat{m}_1(F) = 1 + \sqrt{\eps }D \hat{m}_{1}(\Id) a(G) + \frac{\eps}{2} \hat{Q}( G) + \omega(\sqrt{\eps} G)
\end{align}
with $\sup \left\{\frac{\omega(H)}{|H|^2} : |H|\leq \rho\right\} \rightarrow 0$ as $\rho \rightarrow 0$. 

We concern ourselves with the term $D \hat{m}_{1}( \Id) a(G)$. Recall that $|\hat{m}(R(F)) - \hat{m}(F)| \le C |R(F) - F|^2$ by Lemma \ref{lemma: W,f}(i). For $F = \Id + \sqrt{\eps }G $ this implies by \eqref{eq: Linearization}
\begin{align*}
D \hat{m}_{1}( \Id) a(G) 
&= \e_1 \cdot D\hat{m}(\Id) G = \lim_{\eps \to 0} \e_1 \cdot \frac{\hat{m}(F) - \hat{m}(\Id)}{\sqrt{\eps}} \\
& =  \lim_{\eps \to 0} \e_1 \cdot \frac{\hat{m}(R(F)) - \e_1}{\sqrt{\eps}} + O(\sqrt{\eps}) = \lim_{\eps \to 0} \e_1 \cdot a(G)\e_1 + O(\sqrt{\eps}) = 0.
\end{align*}  
In particular, \eqref{eq: m-lin} then implies $0 \le \frac{1}{\eps} f_\kappa(F) = -\frac{\kappa}{2} \hat{Q}(G) - \frac{1}{\eps}\omega(\sqrt{\eps}G)$ and thus $-\hat{Q}$ is positive semidefinite. We proceed exactly as in the proof of Theorem \ref{th: Gamma-convergence}(i) and conclude 
\begin{align*}
\liminf_{\eps \to 0}\frac{1}{\eps}\int_{\Omega_\eps} \chi_\eps f_{\kappa}(\nabla y'_\eps) &\ge   \liminf_{\eps \to 0} -\int_{\Omega_\eps} \chi_\eps \Big( \frac{\kappa}{2} \hat{Q}( \nabla u'_\eps)  + \frac{\kappa}{\eps} \omega (\sqrt{\eps} \nabla u'_\eps)\Big)
\\ &\ge -\frac{\kappa}{2}\int_{\Omega} \hat{Q}(\nabla u). 
\end{align*}
\eop

\subsection{Recovery sequences}

It remains to construct recovery sequences in order to complete the proof of Theorem \ref{th: Gamma-convergence}. \smallskip

\noindent {\em Proof of Theorem \ref{th: Gamma-convergence}(ii).} 
The basic tool for the proof of the $\Gamma$-limsup-inequality is a density result for $SBV$ functions due to Cortesani and Toader \cite{Cortesani-Toader:1999}. Moreover, a proof very similar to that of Proposition 2.5 in \cite{Giacomini:2005} shows that we may also
impose suitable boundary conditions on the approximating sequence. We suppose ${\cal W}(\Omega, \R^2)$ is the space of all $SBV$ functions $u \in SBV(\Omega, \R^2)$ such that $J_{u}$ is a finite union of (disjoint) segments and $u \in W^{k,\infty}(\Omega \setminus J_{u}, \R^2)$ for all $k$. Then ${\cal W}(\Omega, \R^2)$ is dense in $SBV^2(\Omega, \R^2) \cap L^{\infty}(\Omega, \R^2)$ in the following way: 

For every $u \in SBV^2(\tilde{\Omega}, \R^2) \cap L^{\infty}(\tilde{\Omega}, \R^2)$ with $u=g$ on $\tilde{\Omega} \setminus \Omega$, there exists a sequence $u_n$ and a sequence of neighborhoods $U_n \subset \tilde{\Omega}$ of $\tilde{\Omega}\setminus \Omega$ such that $u_n = g$ on $\Omega_{D,\frac{1}{n}}$ (recall \eqref{eq: Ag}), $u_n \in W^{1,\infty}(U_n)$ and $u_n|_{V_n}  \in {\cal W}(V_n, \R^2)$, where
$V_n \subset \Omega$ is some neighborhood of $\Omega \setminus U_n$, such that $\left\|u_n\right\|_{\infty} \leq \left\|u\right\|_{\infty}$ and
\begin{itemize}
	\item [(i)] $u_{ n} \rightarrow u$ strongly in $L^1(\Omega,\R^2)$, $\nabla u_n \rightarrow \nabla u$ strongly in $L^2(\Omega, \R^2)$, 
	\item[(ii)] $\limsup_{n \rightarrow \infty} \int_{J_{u_{n}}} \phi(\nu_{u_{n}})d{\cal H}^1 \leq \int_{J_{u} }\phi(\nu_{u})d{\cal H}^1$ for every upper semicontinuous function $\phi:S^1 \rightarrow [0,\infty)$ satisfying $\phi(\nu) = \phi(-\nu)$ for every $\nu \in S^1$.
\end{itemize}

Let $u \in SBV^2(\tilde{\Omega}, \R^2)$ with $u=g$ on $\tilde{\Omega} \setminus \Omega$. Without restriction we can assume $u \in  \cap L^{\infty}(\tilde{\Omega}, \R^2)$ as this hypothesis my be dropped by applying a truncation argument and taking $Q(F) \leq C |F|^2$ into account. In fact, it suffices to provide a recovery sequence for an approximation $u_n$ defined above. Although our notion of convergence in Definition \ref{def: convergence} is not given in terms of a specific metric, similarly to a general density result in the theory of $\Gamma$-convergence this can be seen by a diagonal sequence argument. The crucial point is that due to \eqref{eq: D estimate} below we may assume that for $\eps$ sufficiently small (depending on $n$) 
$$ \# {\cal C}^*_{\eps} = \# {\cal D}_{\eps} \leq \frac{C {\cal H}^1(J_{u_n})}{\eps} \le \frac{C {\cal H}^1(J_{u})}{\eps}, $$ 
where $C$ is independent of $n$ and $\eps$. If $(u_{n,\eps})_{\eps}$ is a recovery for $u_n$, one may therefore pass to a diagonal sequence which is a recovery sequence for $u$, in particular converging to $u$ the sense of Definition \ref{def: convergence}. For simplicity write $u$ instead of $u_n$ in what follows.

Let $\delta >0$ and define $J^{\delta}_{u} = \left\{x \in J_{u}, |[u](x)| \geq \delta \right\}$. Since $|[u]|$ is Lipschitz continuous on $J_u$, it cannot oscillate infinitely often between values $\le \delta$ and values $\ge 2\delta$ on a single segment. Consequently, there is a finite number $N_u^{\delta}$ of disjoint subsegments $S_1, \ldots, S_{N_u^{\delta}}$ in $J_u$ such that $|[u]| < 2 \delta$ on every $S_j$ and $|[u]| > \delta$ on $J_u \setminus (S_1 \cup \ldots \cup S_{N_u^{\delta}})$. Note that ${\cal H}^1( \bigcup^{N_u^{\delta}}_{i=1} S_i) \le {\cal H}^1(J_u \setminus J^{2\delta}_u)=: \rho(\delta) \to 0$ for $\delta \to 0$. We cover $S_1, \ldots, S_{N_u^{\delta}}$ by pairwise disjoint  rectangles $Q_1, \ldots Q_{N_u^{\delta}}$ which satisfy $\sum_j {\cal H}^1(\partial Q_i) + |Q_i|  \le C\rho(\delta)$.  It is not hard to see that $|u({ x}) - u(y)| \le C{\cal H}^1(\partial Q_i) + 2\delta$ for $x,y \in Q_j$ as $\nabla u \in L^\infty(\tilde{\Omega})$. 

We modify $u$ on the rectangles $Q_i$: Let $u_\delta = u$ on $\tilde{\Omega} \setminus \bigcup^{N_u^{\delta}}_{i=1} Q_j$ and define $u_\delta =  c_j$ on $Q_j$ for $c_j \in \R^2$ in such a way that $J_{u_\delta} = J^\delta_{u_\delta}$ up to an ${\cal H}^1$-negligible set. As $u \in L^\infty(\tilde{\Omega})$, $\nabla u \in L^\infty(\tilde{\Omega})$ we find $u_\delta \to u$ in $L^1(\tilde{\Omega})$ and $\nabla u_\delta \to \nabla u$ in $L^2(\tilde{\Omega})$. Moreover, we have ${\cal H}^1(J_u \Delta J_{u_\delta}) \le C\rho(\delta) \to 0$ for $\delta \to 0$. 

Consequently, it suffices to establish a recovery sequence for a function $u \in {\cal W}(\Omega)$ with $u = g$ in a neighborhood of $\tilde{\Omega} \setminus \Omega$ and $J_u = J^\delta_u$ for some $\delta >0$. Note after the above modification the segments of $J_u$ might not be pairwise disjoint.

We define $u_{\eps}(x) = u(x)$ for $x \in {\cal L}_{\eps} \cap \tilde{\Omega}$ and let $y_{\eps}(x) = \id + \sqrt{\eps} u_{\eps}(x)$. Clearly we have $u_{\eps} \in {\cal A}_{g_\eps}$ for all $\eps$. By $\tilde{u}_{\eps}, u'_{\eps}$ we again denote the interpolations on $\tilde{\Omega}_{\eps}$. Up to considering a translation of $u$ of order $\eps$, we may assume that $J_{u} \cap {\cal L}_{\eps} = \emptyset$. Let ${\cal D}_{\eps}$ be the sets of triangles where $J_{u}$ crosses at least one side of the triangle. Then 
\begin{align}\label{eq: D estimate}
\# {\cal D}_{\eps} \leq \frac{C {\cal H}^1(J_{u})}{\eps} + C N_u
\end{align} 
for a constant $C > 0$ independent of $u \in {\cal W}(\tilde{\Omega}, \R^2)$ and $\eps$, where $N_u$ denotes the (smallest) number of segments whose union gives $J_{u}$. From now on for the local nature of the arguments we may assume that $J_{u}$ consists of one segment only. Indeed, if $J_u$ consists of segments $S_1,\ldots, S_{N_u}$, which are possibly not disjoint, the number of triangles $\Delta \in \tilde{\cal C}_\eps$ with $\triangle \cap S_{i_1} \cap S_{i_2} \neq \emptyset$ for $1\le i_i < i_2 \le N_u$ scales like $N_u$ and therefore their energy contribution is negligible in the limit. We show 
\begin{equation*}
\bar{{\cal C}}_{\eps} = {\cal D}_{\eps}
\end{equation*}
for $\eps$ small enough. Let $\triangle \in {\cal D}_\eps$. We see that, if $J_{u} = J^{\delta}_{u}$ crosses a spring $\vv$ at point $x_{*}$, say, then a computation similar as in (\ref{eq: [2] triangle}) together with $\nabla u \in L^{\infty}$ shows
\begin{equation}\label{eq: delta jump}
  \left| (\tilde{y}_{\eps})_{\triangle} \, \vv \right| 
  = \left| \frac{1}{\sqrt{\eps}} [u(x_{*})] + O(1) \right| 
  \geq  \frac{\delta}{\sqrt{\eps}} + O(1). 
\end{equation}
Thus, $\triangle \in \bar{{\cal C}}_{\eps}$ for $\eps$ small enough. On the other hand, if we assume $\triangle \notin {\cal D}_{\eps}$, then for at least two springs $\vv \in {\cal V}$ we have $|(\tilde{y}_{\eps})_{\triangle} \, \vv| \leq 1 + \sqrt{\eps} \left\| \nabla u \right\|_{\infty} < 2$ for $\eps$ small enough leading to $\triangle \notin \bar{{\cal C}}_{\eps}$. 

We claim that 
\begin{equation}\label{eq: bondlengths1}
  \Vert\nabla u'_{\eps}\Vert_{L^\infty(\tilde{\Omega})} \le C. 
\end{equation}
This is clear for $\triangle \notin {\cal D}_{\eps} = \bar{\cal C}_{\eps}$ as $\nabla u \in L^{\infty}$. For $\triangle \in  \bar{{\cal C}}_{\eps,3}$ it follows by construction. For $\triangle \in \bar{{\cal C}}_{\eps,2}$ there is a $\vv \in {\cal V}$ such that $(y'_{\eps})_{\triangle} \, \vv = (\tilde{y}_{\eps})_{\triangle} \, \vv = \vv + O(\sqrt{\eps})$. By Lemma \ref{lemma:quadratic lower bound}(i) and (\ref{eq: SO, O}) we get a rotation $R_{\eps} \in SO(2)$ such that
$$ |R_{\eps} - (y'_{\eps})_{\triangle}|^2 
   = \dist^2((y'_{\eps})_{\triangle}, SO(2)) 
   = \dist^2((y'_{\eps})_{\triangle}, O(2)) 
   \leq C W_{\triangle}((y'_{\eps})_{\triangle}) 
   =  O(\eps).$$
This yields $|(y'_{\eps})_{\triangle} - \Id| = O(\sqrt{\eps})$ and thus $|(u'_{\eps})_{\triangle}| = O(1)$.

We note that $\chi_{\tilde{\Omega}_{\eps}} \tilde{u}_{\eps} \to u$ in $L^1$ as $u$ and thus every $\tilde{u}_{\eps}$ is bounded uniformly in $L^{\infty}$ and, $u$ being Lipschitz away from $J_u$, $\tilde{u}_{\eps} \to u$ uniformly on $\tilde{\Omega}_{\eps} \setminus \bigcup_{\triangle \in {\cal D}_{\eps}} \triangle$, where $|\bigcup_{\triangle \in {\cal D}_{\eps}} \triangle| \le C \eps$. Letting ${\cal C}^*_{\eps} = {\cal D}_{\eps}$ this shows that $u_{\eps} \rightarrow u$ in the sense of Definition \ref{def: convergence} recalling (\ref{eq: D estimate}) and the fact that $|(\tilde{u}_{\eps})_{\triangle}| = O(1)$ for $\triangle \notin {\cal D}_{\eps}$. We next establish an even stronger convergence of the derivatives. Consider $\nabla \tilde{u}_{\eps}$ on triangles in ${\cal C}_{\eps} \setminus {\cal D}_{\eps}$. As $u$ is Lipschitz there, the oscillation on such a triangle, $\text{osc}^{\triangle}_{\eps} (\nabla u) := \sup \left\{\left\|\nabla u(x) - \nabla u(x')\right\|_\infty, x,x' \in \triangle\right\}$, tends to zero uniformly (i.e., not depending on the choice of the triangle). We thus obtain 
$$ \int_{\tilde{\Omega}_{\eps} \setminus \cup_{\triangle \in {\cal D}_{\eps}} \triangle} 
   \|\nabla \tilde{u}_{\eps} - \nabla u\|_{\infty}^2 
   \leq \int_{\tilde{\Omega}_{\eps} \setminus \cup_{\triangle \in {\cal D}_{\eps}} \triangle} 
   (\text{osc}^{\triangle}_{\eps} (\nabla u))^2 \rightarrow 0$$
for $\eps \rightarrow 0$, so that even $\chi_{\tilde{\Omega}_{\eps} \setminus \cup_{\triangle \in {\cal D}_{\eps}} \triangle}\nabla \tilde{u}_{\eps} \rightarrow \nabla u$ strongly in $L^2(\tilde{\Omega})$. Note that in fact $\chi_{\tilde{\Omega}_{\eps}} \nabla u'_{\eps} \rightarrow \nabla u$ in $L^2(\Omega)$. Indeed, recall $\# {\cal D}_{\eps} \leq C\eps^{-1}$ by \eqref{eq: D estimate}. Using \eqref{eq: bondlengths1} on the set of broken triangles we then get 
\begin{align*}
\int_{\bigcup_{\triangle \in {{\cal D}}_{\eps}} \triangle} |\nabla u'_{\eps} - \nabla u|^2 \le C \#\bar{{\cal D}}_{\eps} \eps^2 \to 0
\end{align*}
for $\eps \to 0$. 
We now split up the energy in bulk and surface parts 
\begin{align}\label{eq: en spli}
{\cal E}^{\chi}_{\eps}(u_{\eps}) 
      & = {\cal E}^{\rm elastic}_{\eps}(u_{\eps}) + {\cal E}^{\rm crack}_{\eps}(u_{\eps}) + O(\eps) + \frac{1}{\eps}\int_{\Omega_\eps} \chi(\nabla \tilde{y}_\eps) 
\end{align}
as defined in (\ref{eq: auxiliary lower bound}). Note that indeed the contribution $\eps E^{\rm boundary}_{\eps}$ is of order $O(\eps)$ as $\nabla u \in L^\infty(\tilde{\Omega})$ and $J_u \subset \Omega$ since $u = g$ in a neighborhood of $\tilde{\Omega}\setminus \Omega$. We first observe that $\frac{1}{\eps}\int_{\Omega_\eps} \chi(\nabla \tilde{y}_\eps) = 0$ for $\eps$ small enough. Indeed, for $\Delta \in \bar{\cal C}_\eps$ this follows from \eqref{eq: delta jump}. For $\Delta \notin {\cal D}_\eps$ it suffices to recall $|(\tilde{u}_{\eps})_{\triangle}| = O(1)$ which implies that $(\tilde{u}_{\eps})_{\triangle}$ is near $SO(2)$.  Repeating the steps in the elastic energy estimate in (i), applying $\chi_{\Omega_{\eps}} \nabla u'_{\eps} \rightarrow \nabla u$ strongly in $L^2(\Omega)$, \eqref{eq: bondlengths1} and $Q(F) \leq C|F|^2$ for a constant $C > 0$ we conclude that 
\begin{equation}\label{eq: recov1}
 \limsup_{\eps \rightarrow 0}{\cal E}^{\rm elastic}_{\eps}(u_{\eps}) 
   = \frac{4}{\sqrt{3}} \int_{\Omega}  \frac{1}{2} Q(e(u(x))) \, dx. 
   \end{equation}
It is elementary to see that $J_{u}$ crosses 
\begin{equation}\label{eq: spring counting}
{\cal H}^1(J_{u}) \frac{2 |\nu_{u} \cdot \vv|}{\sqrt{3} \eps} + O(1)
\end{equation}
springs in $\vv$-direction for $\vv \in {\cal V}$, where $\nu_u$ is a normal to the segment $J_u$. Recalling (\ref{eq: delta jump}), the crack energy may be estimated by 
\begin{align*}
  &\limsup_{\eps \to 0}{\cal E}^{\rm crack}_{\eps}(u_{\eps}) \\ 
  &\leq \limsup_{\eps \to 0} {\cal H}^1(J_{u})  \ \sup \left\{W(r): r \geq \delta \eps^{-\frac{1}{2}} + O(1) \right\} \frac{2}{\sqrt{3}} \sum_{\vv \in {\cal V}} |\nu_{u} \cdot \vv|  +O(\eps) \\
& = {\cal H}^1(J_{u}) \  \beta  \ \frac{2}{\sqrt{3}} \sum_{\vv \in {\cal V}} |\nu_{u} \cdot \vv|. 
\end{align*}
This together with \eqref{eq: en spli} and \eqref{eq: recov1} shows that $u_{\eps}$ is a recovery sequence for $u$. \eop

Finally, we construct recovery sequences for the functionals ${\cal F}^\chi_\eps$ to conclude the proof of Theorem \ref{th: Gamma-convergence2}.\smallskip

\noindent {\em Proof of Theorem \ref{th: Gamma-convergence2}, second part.} 
Following the proof of Theorem \ref{th: Gamma-convergence}(ii) it suffices to show 
$$\lim_{\eps \to 0}  \frac{1}{\eps}\int_{\Omega_\eps}f_{\kappa}(\nabla \tilde{y}_\eps) = -\frac{\kappa}{2}\int_\Omega \hat{Q}(\nabla  u).$$
First, by \eqref{eq: delta jump} and the definition of $f_{\kappa}$ we get $\int_{\bigcup_{\Delta \in {\cal D}_\eps} \Delta} f_{\kappa}(\nabla \tilde{y}_\eps) = 0$ for $\eps$ small enough. For $\Delta \notin {\cal D}_\eps$ we have $(\nabla \tilde{y}_\eps)_\Delta = (\nabla y'_\eps)_\Delta$ and thus we find $f_{\kappa}((\nabla \tilde{y}_\eps)_\Delta) = - \eps \frac{\kappa}{2} \hat{Q}( ( \nabla u'_\eps)_\Delta) - \kappa\omega(\sqrt{\eps} \nabla (u'_\eps)_\Delta)$ by \eqref{eq: m-lin}. We obtain
\begin{align*}
\frac{1}{\eps}\int_{\Omega_\eps}f_{\kappa}(\nabla \tilde{y}_\eps) &= 
 \frac{1}{\eps}\int_{\Omega_\eps \setminus \bigcup_{\Delta \in {\cal D}_\eps} \Delta} f_{\kappa}(\nabla y'_\eps)  \\
& \le -\frac{\kappa}{2}\int_{\Omega_\eps \setminus \bigcup_{\Delta \in {\cal D}_\eps} \Delta}  \hat{Q}(\nabla u'_\eps)  +  
\frac{C}{\eps}\int_{\Omega_\eps} \omega(\sqrt{\eps} \nabla u'_\eps).
\end{align*} 
Using \eqref{eq: bondlengths1} and the definition of $\omega$ we observe $\frac{1}{\eps}\Vert\omega(\sqrt{\eps} \nabla u'_\eps)\Vert_\infty \to 0$ for $\eps \to 0$. This together with strong convergence $\chi_{\Omega_\eps} \nabla u'_\eps \to \nabla u$ in $L^2(\Omega)$ shows
\begin{align*}
\limsup_{\eps \to 0}  \frac{1}{\eps}\int_{\Omega_\eps}f_{\kappa}(\nabla \tilde{y}_\eps) \le -\frac{\kappa}{2}\int_\Omega \hat{Q}(\nabla u).
\end{align*} 
\eop

\section{Analysis of the limiting variational problem}\label{sec:limiting-variational-problem}

We finally give the proof of Theorem \ref{th: sharpness, uniqueness} determining the minimizers of the limiting functional ${\cal E}$. An analogous result for isotropic energy functionals has been obtained in \cite{Mora:2010}. We thus do not repeat all the steps of the proof provided in \cite{Mora:2010} but rather concentrate on the additional arguments necessary to handle anisotropic surface contributions. 
\smallskip

\noindent {\em Proof of Theorem \ref{th: sharpness, uniqueness}.} 
We first establish a
lower bound for the energy ${\cal E}$. To this end, we begin to estimate $\sum_{\vv \in {\cal V}} |\vv \cdot \nu|$ for $\nu \in S^1$. We recall that $\gamma \in [\frac{\sqrt{3}}{2},1]$ and define $P:[\frac{\sqrt{3}}{2},1] \times S^1 \to [0,\infty)$ by 
\begin{align*}
P(\gamma,\nu) 
  = \begin{cases} \left(1- \sqrt{3} \frac{\sqrt{1 - \gamma^2}}{\gamma}\right) |\vv_\gamma \cdot \nu|, & \gamma >\frac{\sqrt{3}}{2}, \\ \max\big\{ \sqrt{3} |\e_2 \cdot \nu| - |\e_1 \cdot \nu|,0\big\}, & \gamma = \frac{\sqrt{3}}{2}. \end{cases}
 \end{align*}
As $\vv_\gamma$ is unique for $\gamma > \frac{\sqrt{3}}{2}$, the function $P$ is well defined. In the generic case, i.e. for $\gamma > \frac{\sqrt{3}}{2}$, an elementary computation yields
\begin{align*}
  \sum_{\vv \in {\cal V}} |\vv \cdot \nu|
  &\geq |\vv_{\gamma} \cdot \nu| + \sqrt{3} |\vv^{\bot}_{\gamma} \cdot \nu|
   = |\vv_{\gamma} \cdot \nu|
     + \sqrt{3} \left| \pm \frac{1}{\gamma} \e_1 \cdot \nu
     \pm \frac{\sqrt{1 - \gamma^2}}{\gamma} \vv_{\gamma} \cdot \nu
\right| \\
  &\ge \frac{\sqrt{3}}{\gamma}|\e_1 \cdot \nu| + P(\gamma, \nu) 
\end{align*}
for $\nu \in S^1$. In the
first step we used that $\sum_{\vv \in {\cal
V}\setminus\left\{\vv_{\gamma}\right\}} \vv = \pm
\sqrt{3}\vv^{\bot}_{\gamma}$. In the special case $\phi = 0 \Leftrightarrow \gamma = \frac{\sqrt{3}}{2}$, i.e. $\vv_1 = \e_1$, $\vv_{2,3} = \pm \frac{1}{2} \e_1 + \frac{\sqrt{3}}{2} \e_2$ we obtain $\sum_{\vv \in {\cal V}} |\vv \cdot \nu| = |\e_1\cdot \nu| + \sqrt{3} |\e_2 \cdot \nu|$ for $|\nu_2|>\frac{1}{2}$ and $\sum_{\vv \in {\cal V}} |\vv \cdot \nu| = 2|\e_1\cdot \nu|$ for $|\nu_2|\le\frac{1}{2}$, $\nu \in S^1$. Consequently, it is not hard to see that 
$$\sum_{\vv \in {\cal V}} |\vv \cdot \nu| 
  \ge \frac{\sqrt{3}}{\gamma} |\e_1 \cdot \nu| + P(\gamma,\nu)$$
also holds for $\gamma = \frac{\sqrt{3}}{2}$. Thus, we get
$$ {\cal E}(u)
   \geq  \frac{4}{\sqrt{3}} \int_{\Omega}  \frac{1}{2} Q(e(u(x))) \, dx
   + \int_{J_{u}}  \frac{2\beta}{\gamma} |\e_1 \cdot \nu_{u}|
   + \frac{2\beta}{\sqrt{3}} P(\gamma,\nu_u) \, d{\cal H}^1. $$
By Lemma \ref{lemma:W-linearized-properties} we obtain $\min \lbrace Q(F): \e_1^T F\e_1 =r\rbrace = \frac{\alpha}{2} r^2 $. Then using the slicing method (see, e.g., \cite[Section 3.11]{Ambrosio-Fusco-Pallara:2000}) we get
\begin{align}\label{eq:E-uS}
  {\cal E}(u) \geq  \int^{1}_{0}\left( \int^{l}_{0}
  \frac{\alpha}{\sqrt{3}}\left(\e_1^{T}\nabla u(x_1,x_2)\e_1\right)^2 \, dx_1 
  + \frac{2\beta}{\gamma} \#S^{x_2}(u)\right) \,dx_2 + {\cal E}^{\gamma}(u), 
\end{align}
where $\# S^{x_2}$ denotes the number of jumps on a slice $(0,l) \times
\left\{x_2\right\}$ and
$$ {\cal E}^{\gamma}(u)= \int_{J_{u}} \frac{2\beta}{\sqrt{3}
P(\gamma, \nu_u)} \, d{\cal H}^1.$$
In case $\#S^{x_2}(u) \ge 1$, the inner integral in \eqref{eq:E-uS} is obviously bounded from below by $\frac{2 \beta}{\gamma}$. If $\#S^{x_2}(u) =0$, by applyig Jensen's inequality we find that this term is bounded from below by $\alpha l a^2$ due to the boundary conditions. We thus obtain $\inf {\cal E} \geq \min \big\{\frac{\alpha l a^2}{\sqrt{3}}, \frac{2 \beta}{\gamma}\big\}$. On the other hand, it is straighforward to check that ${\cal E}(u^{\rm el}) = \alpha l a^2$ and ${\cal E}(u^{\rm cr}) = \frac{2 \beta}{\gamma}$, which shows that $u^{\rm el}$ is a minimizer for $a < a_{\rm crit}$ and $u^{\rm cr}$ is a minimizer for $a > a_{\rm crit}$. It remains to prove uniqueness: \smallskip

(i) Let $a < a_{\rm crit}$ and $u$ be a minimizer of ${\cal E}$. Since ${\cal E}(u) = {\cal E}(u^{\rm el})$ we infer from \eqref{eq:E-uS} that $u$
has no jump on a.e.\ slice $(0,l) \times \left\{x_2\right\}$ and satisfies $\e_1^{T} \nabla u \, \e_1 = a$ a.e.\ by the imposed boundary values and strict convexity of the
mapping $t \mapsto t^2$ on $[0, \infty)$. Thus, if $J_{u} \neq
\emptyset$,  a crack normal must satisfy $\nu_{u} = \pm \e_2$ ${\cal H}^1$-a.e.
Taking ${\cal E}^{\gamma}(u)$ and the fact that $P(\gamma,\e_2)>0$ for $\gamma \in [\frac{\sqrt{3}}{2},1]$ into account, we then may assume $J_{u} =
\emptyset$ up to an ${\cal H}^1$ negligible set, i.e., $u \in
H^1(\Omega)$.  We find $u_1(x_1,x_2) = a x_1 + f(x_2)$ a.e.\ for a
suitable function $f$, and the boundary condition $u_1(0,x_2)=0$ yields $f = 0$ a.e. In particular, $\e_1^{T}
\nabla u \, \e_2 = 0$ a.e. Applying strict convexity of $Q$ on symmetric
matrices (Lemma \ref{lemma:W-linearized-properties}) we now observe
$\e_2^{T} \nabla u \, \e_2 = -\frac{a}{3}$ and $\e_1^{T} \nabla u \, \e_2
+ \e_2^{T} \nabla u \, \e_1 = 0$ a.e. So the derivative has the form
$$ \nabla u(x) = \begin{footnotesize} \begin{pmatrix} a & 0  \\ 0 &
-\frac{a}{3} \end{pmatrix} \end{footnotesize} \text{ for a.e.\ $x$}. $$
Since $\Omega$ is connected, we conclude $u(x)= (0,s) + F^{a} x = u^{\rm
el}(x)$ a.e.

(ii) Let $a > a_{\rm crit}$, $\phi \neq 0$ and $u$ be a minimizer of ${\cal E}$. 
We again consider the lower bound  \eqref{eq:E-uS} for the energy ${\cal E}$ and now obtain that on
a.e.\ slice $(0,l) \times \left\{x_2\right\}$ a minimizer $u$ has  precisely one jump
and that $\e_1^{T} \nabla u \, \e_1 = 0$ a.e. Now Lemma \ref{lemma:W-linearized-properties} shows that $\nabla u$ is antisymmetric a.e. As a consequence, the linearized rigidity estimate for SBD functions of Chambolle, Giacomini and Ponsiglione \cite{Chambolle-Giacomini-Ponsiglione:2007} yields that there is a Caccioppoli partition $(E_i)$ of $\Omega$ such that 
$$ u(x) = \sum_{i} (A_i x + b_i) \chi_{E_i} \quad \text{and} \quad 
   J_u = \bigcup_{i} \partial^* E_i, $$
where $A_i^T = -A_i \in \R^{2 \times 2}$ and $b_i \in \R^2$. (See \cite{Ambrosio-Fusco-Pallara:2000} for the definition and basic properties of Caccioppoli partitions.) As ${\cal E}^{\gamma}(u) = 0$, we also note that $\nu_u \perp \vv_{\gamma}$ a.e.\ on $J_u$. Following the arguments in \cite{Mora:2010}, in particular using regularity results for boundary curves of sets of finite perimeter and exhausting the sets $\partial^* E_i$ with Jordan curves, we find that 
$$ J_u = \bigcup_{i} \partial^* E_i  \subset (p, 0) + \R \vv_{\gamma} $$ 
for some $p$ such that $(p, 0) + \R \vv_{\gamma}$ intersects both segments $(0, l) \times \{0\}$ and $(0, l) \times \{1\}$. We thus obtain that $(E_i)$ consists of only two sets: $E_1$ to the left and $E_2$ to the right of $(p, 0) + \R \vv_{\gamma}$, say. Due to the boundary conditions we conclude that $A_1 = A_2 = 0$ and $b_1 = (0,s)$, $b_2 = (al, t)$ for suitable $s, t \in \R$. 

(iii) Let $a > a_{\rm crit}$, $\phi = 0$ and $u$ be a minimizer of ${\cal E}$. We follow the lines of the proof in (ii). The only difference is that ${\cal E}^\gamma(u) = 0$ now implies that $|\nu_u \cdot \e_1| \ge \frac{\sqrt{3}}{2}$ a.e. and then arguing similarly as before we obtain 
$$J_u \subset h((0,1))$$
up to an ${\cal H}^1$-negligible set, where $h:(0,1) \to [0,l]$ is a Lipschitz function with $|h'| \le \frac{1}{\sqrt{3}}$ a.e. We now conclude as in (ii). 
\eop


 \typeout{References}


\begin{thebibliography}{10}





\bibitem{Alberti-Mantegazza}
{\sc G.~Alberti, C.~Mantegazza}. 
\newblock {\em A note on the theory of $SBV$ functions}.
\newblock Boll.\ Un.\ Mat.\ Ital.\ B(7)\ 
\newblock {\bf 11} (1989), 375-–382. 


\bibitem{Alicandro-Focardi-Gelli:2000}
{\sc R.~Alicandro, M.~Focardi, M.~S.~Gelli}. 
\newblock {\em Finite-difference approximation of energies in fracture mechanics}. 
\newblock Ann.\ Scuola Norm.\ Sup.\ 
\newblock {\bf 29} (2000), 671--709. 

\bibitem{Ambrosio:89}
{\sc L.~Ambrosio}. 
\newblock {\em A Compactness Theorem for a Special Class of Functions of Bounded
Variation}.
\newblock Boll.\ Un.\ Mat.\ Ital.\ 
\newblock {\bf 3-B} (1989), 857-–881. 

\bibitem{Ambrosio:90}
{\sc L.~Ambrosio}. 
\newblock {\em Existence theory for a new class of variational problems}.
\newblock Arch.\ Ration.\ Mech.\ Anal.\
\newblock {\bf 111} (1990), 291--322. 

\bibitem{Ambrosio-Fusco-Pallara:2000} 
{\sc L.~Ambrosio, N.~Fusco, D.~Pallara}.
\newblock {\em Functions of bounded variation and free discontinuity problems}. 
\newblock Oxford University Press, Oxford 2000. 





\bibitem{Braides:02}
{\sc A.~Braides}.
\newblock {\em $\Gamma$-convergence for Beginners}.
\newblock Oxford University Press, Oxford 2002.

\bibitem{Braides-Cicalese:2007} 
{\sc A.~Braides, M.~Cicalese}. 
\newblock {\em Surface energies in nonconvex discrete systems}.
\newblock Math.\ Models Methods Appl.\ Sci.\
\newblock {\bf 17} (2007), 985--1037. 

\bibitem{Braides-DalMaso-Garroni:1999} 
{\sc A.~Braides, G.~Dal Maso, A.~Garroni}. 
\newblock {\em Variational formulation of softening phenomena in fracture mechanics. The one-dimensional case}.
\newblock Arch.\ Ration.\ Mech.\ Anal.\
\newblock {\bf 146} (1999), 23--58. 

\bibitem{Braides-Gelli:2002-1} 
{\sc A.~Braides, M.~S.~Gelli}. 
\newblock {\em Limits of discrete systems without convexity hypotheses}. 
\newblock Math.\ Mech.\ Solids
\newblock {\bf 7} (2002), 41--66. 


\bibitem{Braides-Gelli:2002-2} 
{\sc A.~Braides, M.~S.~Gelli}. 
\newblock {\em Limits of discrete systems with long-range interactions}. 
\newblock J.\ Convex Anal.\
\newblock {\bf 9} (2002), 363--399. 



\bibitem{Braides-Lew-Ortiz:06}
{\sc A.~Braides, A.~Lew, M.~Ortiz}.
\newblock {\em Effective cohesive behavior of layers of interatomic planes}.
\newblock Arch.\ Ration.\ Mech.\ Anal.\ 
\newblock {\bf 180} (2006), 151--182.

\bibitem{Braides-Solci-Vitali:07}
{\sc A.~Braides, M.~Solci, E.~Vitali}.
\newblock {\em A derivation of linear elastic energies from pair-interaction atomistic systems}.
\newblock Netw.\ Heterog.\ Media 
\newblock {\bf 2} (2007), 551--567.


\bibitem{Chambolle-Giacomini-Ponsiglione:2007}
{\sc A.~Chambolle, A.~Giacomini, M.~Ponsiglione}. 
\newblock {\em Piecewise rigidity}.
\newblock J.\ Funct.\ Anal.\ Solids 
\newblock {\bf 244} (2007), 134--153. 


\bibitem{Conti-Dolzmann-Kirchheim-Mueller:06}
{\sc S.~Conti, G.~Dolzmann, B.~Kirchheim and S.~M{\"u}ller}. 
\newblock {\em Sufficient conditions for the validity of the Cauchy-Born rule close to
$SO(n)$}.
\newblock J.\ Eur Math.\ Soc.\ (JEMS)
\newblock {\bf 8} (2006), 515--539. 


\bibitem{Cortesani-Toader:1999}
{\sc G.~Cortesani, R.~Toader}. 
\newblock {\em A density result in SBV with respect to non-isotropic energies}.
\newblock Nonlinear Analysis
\newblock {\bf 38} (1999), 585--604. 



\bibitem{DalMaso:93}
{\sc G. Dal Maso}.
\newblock {\em An introduction to $\Gamma$-convergence}.
\newblock Birkh{\"a}user, Boston $\cdot$ Basel $\cdot$ Berlin 1993. 











\bibitem{DeGiorgi-Ambrosio:1988}
{\sc E.~De Giorgi, L.~Ambrosio}. 
\newblock {\em Un nuovo funzionale del calcolo delle variazioni}. 
\newblock Acc.\ Naz.\ Lincei, Rend.\ Cl.\ Sci.\ Fis.\ Mat.\ Natur.\ 
\newblock {\bf 82} (1988), 199--210. 


\bibitem{Federer:1969}
{\sc H.~Federer}.
\newblock {\em Geometric measure theory}. 
\newblock Springer, New York, 1969. 




\bibitem{Focardi-Gelli:2003}
{\sc M.~Focardi, M.~S.~Gelli}. 
\newblock {\em Approximation results by difference schemes of fracture energies: the vectorial case}. 
\newblock NoDEA\ Nonlinear\ Differential\ Equations\ Appl.\ Vol. 4\ No.
\newblock {\bf 4} (2003), 469--495. 




\bibitem{Giacomini:2005}
{\sc A.~Giacomini}. 
\newblock {\em Ambrosio-Tortorelli approximation of quasi-static evolution of brittle fractures}.
\newblock Calc.\ Var.\ Partial\ Differential\ Equations.
\newblock {\bf 22} (2005), 129–-172. 



\bibitem{Francfort-Marigo:1998}
{\sc G.~A.~Francfort, J,~J.~Marigo}. 
\newblock {\em Revisiting brittle fracture as an energy minimization problem}.
\newblock J.\ Mech.\ Phys.\ Solids 
\newblock {\bf 46} (1998), 1319--1342. 


\bibitem{FriedrichSchmidt:2014}
{\sc M.~Friedrich, B.~Schmidt}.
\newblock {\em An atomistic-to-continuum analysis of crystal cleavage in a two-dimensional model problem}. 
\newblock  J.\ Nonlin.\ Sci.\ 
\newblock {\bf 24} (2014), 145--183. 

\bibitem{FrieseckeJamesMueller:02}
{\sc G.~Friesecke, R.~D.~James, S.~M{\"u}ller}.
\newblock {\em A theorem on geometric rigidity and the derivation of nonlinear plate theory from three-dimensional elasticity}. 
\newblock Comm.\ Pure Appl.\ Math.\ 
\newblock {\bf 55} (2002), 1461--1506. 


\bibitem{Friesecke-Theil:02}
{\sc G.~Friesecke, F.~Theil}. 
\newblock {\em Validity and failure of the Cauchy-Born hy\-po\-the\-sis in a two-dimensional mass-spring lattice}. 
\newblock J.\ Nonlinear Sci.\ 
\newblock {\bf 12} (2002), 445--478.


\bibitem{Mora:2010}
{\sc C.~Mora-Corral}. 
\newblock {\em Explicit energy-minimizers of incompressible elastic brittle bars under uniaxial extension}.
\newblock C.\ R.\ Acad.\ Sci.\ Paris 
\newblock {\bf 348} (2010), 1045--1048. 


\bibitem{Negri:2003}
{\sc M.~Negri}. 
\newblock {\em Finite element approximation of the {G}riffith's model in fracture mechanics}.
\newblock Numer.\ Math.\ 
\newblock {\bf 95} (2003), 653--687. 




\bibitem{Schmidt:2009}
{\sc B.~Schmidt}. 
\newblock {\em On the derivation of linear elasticity from atomistic models}. 
\newblock Netw.\ Heterog.\ Media 
\newblock {\bf 4} (2009), 789--812. 















\end{thebibliography}
\end{document}